\documentclass[11pt]{amsart}
\usepackage{amssymb}
\usepackage{graphicx}
\usepackage{color}
\usepackage{here}
\input xypic
\xyoption{all}
\newcommand{\BR}{{\mathbb{R}}}
\newcommand{\BC}{{\mathbb{C}}}
\newcommand{\BF}{{\mathbb{F}}}
\newcommand{\BP}{{\mathbb{P}}}
\newcommand{\gD}{\Delta}
\newcommand{\gd}{\delta}
\newcommand{\gb}{\beta}
\newcommand{\gc}{\gamma}
\newcommand{\gC}{\Gamma}
\newcommand{\gs}{\sigma}
\newcommand{\gS}{\Sigma}
\newcommand{\gO}{\Omega}
\newcommand{\gom}{\omega}
\newcommand{\gl}{\lambda}
\newcommand{\ga}{\alpha}

\newcommand{\gth}{\theta}
\newcommand{\gTh}{\Theta}
\newcommand{\cM}{{\mathcal{M}}}
\newcommand{\cO}{{\mathcal{O}}}
\newcommand{\caR}{{\mathcal{R}}}
\newcommand{\caL}{{\mathcal{L}}}

\newcommand{\ti}[1]{\tilde{#1}}
\newcommand{\ol}[1]{\overline{#1}}
\newcommand{\Pic}{\mathrm{Pic}}
\newcommand{\jac}{\mathrm{Jac}}
\newcommand{\prym}{\mathrm{Prym}}
\newcommand{\gal}{\mathrm{Gal}}
\newcommand{\Sym}{\mathrm{Sym}}
\newcommand{\odd}{\mathrm{odd}}
\newcommand{\even}{\mathrm{even}}
\newcommand{\half}{\frac{1}{2}}
\newcommand{\sm}{\smallsetminus}
\newcommand{\hra}{\hookrightarrow}
\newcommand{\ra}{\longrightarrow}

\swapnumbers
\theoremstyle{plain}
\newtheorem{lemma}{Lemma}[section]
\newtheorem{theorem}[lemma]{Theorem}
\newtheorem{proposition}[lemma]{Proposition}
\newtheorem{corollary}[lemma]{Corollary}
\theoremstyle{definition}
\newtheorem*{remark}{Remark}
\newtheorem{prd}[lemma]{Proposition-Definition}

\newtheorem{ntt}[lemma]{Notations}

\newtheorem{dsc}[lemma]{}

\newcommand{\href}[2]{#2}
\begin{document}

\title[The genus 3 AGM]{An explicit formula for the arithmetic geometric mean in genus $3$}
\author{D. Lehavi}
\thanks{D. Lehavi was partially
supported by Israel-US BSF grant 1998265.}
\address{Mathematics Department, Princeton University, Fine Hall,
Washington Road, Princeton, NJ 08544, USA}
\email{dlehavi@math.princeton.edu}
\author{C. Ritzenthaler}
\thanks{C. Ritzenthaler
acknowledges the financial support provided through the European
Community's Human Potential Programme under contract HPRN-CT-2000-00114,
GTEM}
\address{Institut de Math\'ematiques de Luminy\\
163 Avenue de Luminy, Case 907 \\
13288 Marseille, France}
\email{ritzenth@math.jussieu.fr}
\date{\today}
\keywords{Prym varieties, Arithmetic Geometric Mean}
\subjclass{14H40,14H45,14Q05}
\begin{abstract}
The arithmetic geometric mean algorithm for calculation of elliptic integrals
of the first type
was introduced by Gauss. The analog algorithm for Abelian integrals
of genus $2$ was introduced by Richelot (1837) and Humbert (1901). We present
the analogous algorithm for Abelian integrals of genus $3$.
\end{abstract}
\maketitle

%
\section{Introduction}\label{Sintro}
%
\begin{dsc}
The Arithmetic Geometric Mean (AGM) was discovered by Lagrange in 1785 and
independently by Gauss in 1791. It is described as follows:
given two positive numbers $a$ and $b$, define $M(a,b)$ as the limit of the following convergent sequences:
\[
  \begin{array}{ll}
     a_0:=a,                    &b_0=b \\
     a_{n+1}=\frac{a_n+b_n}{2},  &b_{n+1}=\sqrt{a_n b_n}.
  \end{array}
\]
During the period 1791-1799 Gauss discovered a relation between the AGM and
elliptic curves:
\end{dsc}
\begin{theorem}[Gauss, see \cite{Cox},\cite{BM}]\label{TGagM}
For each pair of positive real numbers in the AGM double sequence $a_n>b_n>0$
define:
\[
  e_{n1}:=\frac{1}{3}(a_n^2+b_n^2),\quad e_{n2}:=e_{n1}-b_n^2,\quad
  e_{n3}:=e_{n1}-a_n^2.
\]
Denote by $E_n$ the elliptic curve given by the equation
\[
  y_n^2=4(x_n-e_{n1})(x_n-e_{n2})(x_n-e_{n3}),
\]
then the following sequence of Abelian groups:
\[
0\rightarrow\{0,(e_{n1},0)-\infty\}\hra E_n
\stackrel{\frac{dx_n}{y_n}\mapsto\frac{dx_{n+1}}{y_{n+1}}}
{-\!\!-\!\!\!\ra} E_{n+1}\rightarrow0
\]
is exact for all $n$.
\end{theorem}
\begin{dsc}
It is easy to see that we create in this way a sequence of
$2$-isogenous elliptic curves.
Gauss generalized the definition of the AGM to the complex numbers;
in this case  there
is a choice involved when taking the square root. Gauss described the
resulting correspondence, see \cite{Cox} for the description.\\
Recall that the real points in the Picard group of a curve which is
defined over $\BR$ are the
divisor classes which are invariant under the action of the group
$\gal(\BC/\BR)$ (see \cite{GH} sections 1-5).
Gauss proved that if the $2$-torsion points of $\Pic(E)$ are real
then there is a unique $\ga\in\Pic(E)[2]$ such that the $2$-torsion points
of $\Pic(E'_\ga)[2]$ are real.
Applying this property iteratively, one gets an algorithm for calculating
elliptic integrals (see \cite{Cox},\cite{BM}) of the form:
\[
  \int^{e_2}_{e_3}\frac{dx}{\sqrt{(x-e_{1})(x-e_{2})(x-e_{3})}}.
\]
This iterative algorithm is applied in numeric evaluations of certain
types of Abelian integrals (see \cite{BB}).

In genera higher than $1$ one can hope for an isogeny between
the Jacobians of the curves. Before stating results in higher
genera, we recall some
facts on polarized Abelian varieties.
\end{dsc}
\begin{dsc}\label{DPPAV}
A pair $(A,\gTh)$ where $A$ is an Abelian variety and $\gTh$ is a
divisor of $A$ is called
a {\em polarized} Abelian variety. The divisor $\gTh$ is called the
{\em theta divisor} of the polarized Abelian variety and the map
$A\ra\Pic^0(A)$ defined by $a\mapsto T_a^{-1}(\gTh)-\gTh$ is
called the {\em polarization} of the pair $(A,\gTh)$.
Since the kernel of the polarization is a finite Abelian group
whose rank is bounded by twice the genus of $A$ we describe its isomorphism type by
$2$ copies of a monotonic sequence of $\mathrm{genus}(A)$ natural numbers. We
will abuse of notations and denote  for instance the sequence
$(2,2,2,2,1,1),(2,2,2,2,1,1)$ by $2^4 1^2$. If the
polarization type is trivial we say that the polarization is {\em principal};
in this case we say that $A$ is a principally polarized Abelian or PPAV.
The translates of $\gTh$ which contain $0$ are called the
{\em theta characteristics} of
$A$. The theta characteristics are called even or odd if their multiplicity at
$0$ is even or odd.
The theta characteristics of an Abelian variety $A$ induce a symplectic
structure on the $\BF_2$-vector space $A[2]$
(the group of $2$-torsion points in $A$) in the following
way: Let $\gth$ be a theta characteristic of $A$, then the map
\[
  \begin{aligned}
  q:A[2]&\ra     \BF_2 \\
       a&\mapsto h^0(\gth+a)+h^0(a)\mod 2
  \end{aligned}
\]
is a quadratic form over $\BF_2$. The quadratic form $q$ induces
the following symplectic pairing on the group $A[2]$:
\[
  \langle a,b\rangle=q(a+b)-q(a)-q(b)\mod 2.
\]
This pairing is called the {\em Weil pairing}. If $G$ is a subgroup of $A[2]$
we will use the notation
\[
  G^\perp:=\{a\in A[2]|\langle a,g\rangle=0 \text{ for all }g\in G\}.
\]
If $A$ is a PPAV and $G$ is a subgroup of the group $A[2]$, then $\Pic^0(A/G)$
is naturally isomorphic to $\Pic^0(A)/G^\perp$. Whence the Abelian
variety $A/G$ is principally polarized only if $G$ is $0,A[2]$, or
a maximal isotropic group of $A[2]$ with respect to the Weil pairing.
Finally recall that a Jacobian of a smooth curve $C$ is principally polarized
by the theta divisor
$\gTh_C$ - the image of $\Sym^{g-1}C$ in $\jac(C)$ under the Abel map.
\end{dsc}
\begin{dsc}
The dimension of the moduli space of principally polarized Abelian varieties of genus $g$ is $g(g+1)/2$, while
the dimension of the moduli space of curves of genus $g \geq 2$ is $3g-3$. Thus if $C$ is a generic
curve of genus $g\leq 3$ and $L$ is a maximal isotropic subgroup of $\jac(C)[2]$ then
$\jac(C)/L$ is a Jacobian of some curve $C'$. This motivates the following
questions:
\begin{itemize}
\item Is this an algebraic correspondence ?
\item Does there exist a curve $C'$ for {\em every} pair $C,L$ in genera $2,3$ (and not only generically)?
\item If $C$ is a real curve, are there ``distinguished'' maximal subgroups ?
\item What is the situation for genera higher then $3$ ?
\end{itemize}
In the case $g=2$ the first three  questions were settled by Richelot (see
\cite{Ric}) and Humbert (see
\cite{Hum}). Both  Richelot \& Humbert described an iterative algorithm for
curves of
genus $2$ with $6$ real Weierstrass points, Richelot described
{\em algebraically} the curve
$C'$ in terms of the curve $C$ and Humbert described
the isomorphism
\[
  H^0(C,K_C)\ra H^0(C',K_{C'}).
\]
See \cite{DL} section 4 for a modern review of the construction and \cite{BM}
for the resulting integration
identities and the iterative integration algorithm resulting in the real case.
Using modern tools (namely GAGA), the answer to the first question is
immediately positive in all genera.
Donagi \& Livn\'e solved the second question for genus $g=3$,
and answered negatively the last question:
\end{dsc}
\begin{theorem}[Donagi \& Livn\'e, see \cite{DL}]\label{TDL}
Let $C$ be a smooth curve of genus $g$ over a base field of
characteristic different from $2,3$. Let $L\subset\jac(C)[2]$ be a
maximal isotropic subgroup (with respect to the Weil pairing) then:
\begin{itemize}
\item If $g=3$ then there exists a curve $C'$ such that
$\jac(C')\cong\jac(C)/L$. The curve $C'$ can be described algebraically in
terms of the curve
$C$ and the maximal isotropic subgroup $L$.
\item If $g>3$ then generically there is {\em no} curve $C'$ such that
$\jac(C')\cong\jac(C)/L$.
\end{itemize}
\end{theorem}
\begin{dsc}
The proof that Donagi and Livn\'e presented for genus $3$ is
constructive in the set-theoretic sense. However, as a basis for
explicit work it has drawbacks; it is not clear how to give
coordinates to the spaces and functions involved or how one can
track the canonical classes.

The object of this paper is to extend Gauss's original work on
curves of genus $1$ to the case of genus $3$. In Section
\ref{Sconst} we use the Coble-Recillas construction (reviewed in
\ref{Ctrigo}-\ref{TPC}) to give an alternative construction to the
one proposed in \cite{DL}; we describe the curve $C'$ in terms of
the pair $(C,L)$, where the curve $C$ is a generic  curve. In
Section \ref{Scanon} we describe the isomorphism between the
canonical classes $|K_C|$ and $|K_{C'}|$. In Section \ref{Sprog} we
derive the formulas describing the curve $C'$ in terms of the curve
$C$ and the isomorphism $H^0(C,K_C)\cong H^0(C',K_{C'})$. We also
show how to iterate the construction. In Section \ref{Sreal} we
concentrate on real curves: assuming the curve $C$ is a real M-curve
(i.e. a real curve with $4$ components) we present a distinguished
maximal isotropic subgroup of
$\jac(C)[2]$, resulting in an iterative integration algorithm.

Although our construction is stated over the complex numbers it is
mostly algebraic. The complex structure is used in one crucial
point: we use a result (due to Jordan in \cite{Jo} section 332, or
see Harris's modern approach in \cite{Ha}) stating that the Galois
group of the bitangents of a smooth plane quartic over the complex
numbers is $SP_6(2)$. In the rest of the paper  we require only that
the characteristic of the base field is the one arising from using
the bigonal and trigonal constructions; the characteristic of the
base field $K$ is not $2$ or $3$ (see the discussion in the
introduction to \cite{DL}).

Some of the proofs in Sections \ref{Sprog},\ref{Sreal} are computer
aided proofs. The {\em Mathematica} and {\em MAGMA}  programs which
generated the computer part of the proofs appears in \cite{LR}.
\end{dsc}
%
%
\begin{dsc}\label{Dfinite}
We remark that the AGM has a nice
application in the area of curves over finite fields:
Mestre observed that the theta-function identities involved in the AGM over
$p$-adic fields can be used to study the number of points in curves over
finite fields.
See \cite{Me},\cite{LL} for results elliptic and hyperelliptic
results and \cite{Rit} for results
on non hyperelliptic genus $3$ curves.
\end{dsc}
\subsection*{Acknowledgments}
Section \ref{Sconst} and \ref{Scanon} of this paper were part of the first
author's Ph.D. thesis, done under the advising of R. Livn\'e, who suggested
the series of questions above. I. Dolgachev kindly gave the first author
an early copy of his preprint \cite{Dol2}, and introduced him to the works of
A. Coble.
%
\section{The construction}\label{Sconst}
%
\begin{dsc} \label{flagLp}
The idea behind our construction is to filter the level-$2$ data and
perform the construction in three steps, using {\em non} principally
polarized Abelian varieties to keep track of the level data.
Throughout this paper we fix a {\em generic} curve $C$ of genus $3$,
and a maximal isotropic flag $\caL=(\langle \ga
\rangle=L_1\subsetneq L_2\subsetneq L_3)$ with respect to the Weil
pairing on $\jac(C)[2]$. Let $C'$ be a curve such that
$\jac(C')\cong\jac(C)/L_3$. The flag $\caL$ induces a dual isotropic
flag $\caL'=(L'_1\subsetneq L'_2\subsetneq L'_3)$ in $\jac(C')[2]$
in the following way: $L_1'$ (resp. $L_2'$, resp. $L_3'$) is the
image of $L_2^{\perp}$ (resp. $L_1^{\perp}$, resp. $\jac(C)[2]$)
under the map $\jac(C) \to \jac(C)/L_3$. We denote by $\ga$
(respectively $\ga'$) the non trivial element in $L_1$ (respectively
$L'_1$). Using the Coble-Recillas construction we will introduce
ramified double covers $Y\ra E$ and $Y'\ra E'$ such that there are
natural isomorphisms
\[
  \prym(Y/E)\cong\jac(C)/L_1^\perp,\quad
   \prym(Y'/E')\cong\jac(C')/{L'_1}^\perp.
\]
Using a bigonal construction (see \cite{Don} p. 68-69 for an overview of
the bigonal construction) we will prove that the polarized Abelian varieties
$\prym(Y'/E')$ and $\prym(Y/E)$ are dual to one another, up to finite
data arising from $L_2$.

\end{dsc}
\begin{dsc}
Before describing the geometry of the construction (in
\ref{TPC}-\ref{TAGM}), we describe the finite symplectic algebra
involved - the level-2 structure of the curve $C$. Define the
following $SP_6(2)$ equivariant surjective map:
\[
  \begin{aligned}
  D:\{\substack{\text{unordered pairs of distinct}\\
                \text{odd theta characteristics of }C}\}
  &\ra\jac(C)[2]\sm\{0\} \\
  \{\theta_1,\theta_2\}&\mapsto \theta_1-\theta_2,
  \end{aligned}
\]
(recall that the group $SP_6(2)$ acts on the odd theta
characteristics via the monodromy action, see \cite{Ha}). Since the
group $SP_6(2)$ acts transitively on the set $\jac(C)[2]\sm\{0\}$,
all the fibers of the map $D$ are of the same order:
$\binom{28}{2}/63=6$. For any $2$-torsion point $\gc$ in $\jac(C)$
we define the {\em Steiner system} $$\gS_\gc:=\cup_{x \in
D^{-1}(\gc)} \{x\}.$$ Since the pairs in a fiber of the map $D$ do
not intersect, the order of all the Steiner systems is $12$. In
Propositions \ref{P46} and \ref{Plev} below we discuss the relation
between the symplectic structure on the vector space $\jac(C)[2]$
and the combinatorics of the $63$ Steiner systems.
\end{dsc}
\begin{proposition}\label{P46}
Let $\gc,\gc'$ be two distinct elements in $\jac(C)[2]\sm\{0\}$,
then $\#(\gS_\gc\cap\gS_{\gc'})$ is $4$ if
$\langle\gc,\gc'\rangle=0$, and $6$ otherwise.
\end{proposition}
\begin{proof}
Since the group $SP_6(2)$ acts transitively on pairs of distinct
elements in $\jac(C)[2]\sm\{0\}$ with the same Weil pairing, the
order of the set $\gS_\gc\cap\gS_{\gc'}$ depends only on the Weil
pairing $\langle \gc,\gc'\rangle$. We denote the two possible
intersection orders by $n_0,n_1$. Any odd theta characteristic sits on
$28-1=27$ different $\gS_\gc$s and the number of $\gc\in\jac(C)[2]\sm\{0,\ga\}$
such that $\langle\gc,\ga\rangle=0$ (resp. $1$) is $30$ (resp. $32$). So we get
\[
  \begin{aligned}
  12 \cdot 27
 =& \sum_{\ga\ne 0} \#\{\gth \;|\; \gth\in \gS_{\ga} \cap \gS_{\gc} \}  \\
 =& \# \gS_{\ga} + \#\{\ga \; |\; \langle\gc,\ga \rangle=0\} n_0
  +\#\{\ga \; |\; \langle\gc,\ga \rangle=0\} n_1\\
=&12+30n_0+32n_1,
  \end{aligned}
\]
and the unique non-negative integer solution of this equation is $n_1=6,n_0=4$.
\end{proof}
\begin{prd}\label{Plev}
The following properties hold:
\item $\bullet$
Let $\gc$ be a $2$-torsion point in $\jac(C)$ then the map
\[
  \begin{aligned}
  D_\gc:\{\substack{\text{unordered pairs of distinct}\\
                \text{classes in }\gS_\gc/\gc}\}&\ra\jac(C)[2]/\gc\\
   \{a_i,a_j\}   &\mapsto a_i-a_j.
  \end{aligned}
\]
is an isomorphism on $(\gc^\perp/\gc)\sm\{0\}$. Moreover, the map
$D_\gc$ maps the intersection pairing to the Weil pairing.
\item $\bullet$
Let $H$ be an isotropic subgroup of $\jac(C)[2]$ of order $4$, then
there is an unique theta characteristic  $\theta$  such that for all
$\ga \in H$, $\theta + \ga$ is an odd theta characteristic.
Denoting the set $\{\theta+\ga, \ga \in H\}$ by $\gC_H$, one
has $\gC_H=\cap_{\gc\in H\sm\{0\}}\gS_\gc$.
\item $\bullet$
A maximal isotropic subgroup of $\jac(C)[2]$ containing $\ga$ is
represented by a partition of $\gS_{\ga}/\ga$ in three
pairs.
\item $\bullet$
If $7$ Steiner systems intersect at $7$ mutually distinct $4$-tuples,
then there is a maximal isotropic subgroup $G\subset\jac(C)[2]$
such that these $7$ Steiner systems are the Steiner systems of the non-zero
elements of $G$.
\end{prd}
\begin{proof}
Let $G_\gc\subset SP_6(2)$ be the stabilizer of $\gc$,
then the $G_\gc$ orbits of $\jac(C)[2]/\gc$ are the sets
\[
  \{0\},(\gc^\perp/\gc)\sm\{0\},(\jac(C)[2]\sm\gc^\perp)/\gc,
\]
which are of orders $1,15,16$ respectively. Since the map $D_\gc$ is
a non-trivial $G_\gc$ equivariant map, and since the order of the
set of unordered pairs of distinct points in $\gS_\gc$ is
$\binom{6}{2}=15$, the map $D_\gc$ is a 1-1 map on the projective
space $(\gc^\perp/\gc)\sm\{0\}$. By a similar counting argument we
prove the claim on the Weil pairing. The last three assertions
follow from Proposition \ref{P46}.
\end{proof}
\begin{dsc}\label{Dbasis}
Using the description of {\em one} maximal isotropic flag in terms
of odd theta characteristics we describe the combinatorics of {\em
two} isotropic flags: $\caL'=(\langle \ga' \rangle=L_1'
\subsetneq L_2'\subsetneq L_3')$ (see \ref{flagLp}) and
$\ti{\caL}=(\langle \ti{\ga}
\rangle=\ti{L}_1\subset\ti{L}_2\subset\ti{L_3})$ in $\jac(C')[2]$
such that:
\[
  \ti{L}_1\oplus L'_3={L'}_2^\perp,\quad \ti{L}_2\oplus L'_3={L'}_1^\perp,\quad \ti{L}_3\oplus L'_3=\jac(C')[2].
\]
Such a description is essential for iterating the algorithm, as the
pair $(C',\ti{\caL})$ should be the starting point of the second
iteration, playing the same role that the pair $(C,\caL)$ played in
the first iteration. Note that $\ti{\caL}$ is not uniquely defined.
For further uses, we will need the following lemma:
\end{dsc}
\begin{lemma}\label{L2t}
For any subgroup $H \subset \ti{L}_3$ of order $4$ one has $\# \gC_H
\cap\gS_{\ga'}=0$ if $H=\ti{L}_2$ and $\# \gC_H\cap\gS_{\ga'}=2$
otherwise.
\end{lemma}
\begin{proof}
Since we have $\ti{L}_2\subset{L'}_1^\perp$ the Weil pairings of the
non trivial elements in
$\ti{L}_2/\ti{L}_1,(L'_1\oplus\ti{L}_1)/\ti{L}_1$ are all $0$. By
Proposition \ref{Plev} the intersection
\[
D_{\ti{\ga}}^{-1}(\ti{L}_2/\ti{L}_1)\cap
  D_{\ti{\ga}}^{-1}((L'_1\oplus\ti{L}_1)/\ti{L}_1)
\]
is empty. Whence, The intersection
\[
 \gC_{\ti{L}_2}\cap\gS_{\ga'}=\gC_{\ti{L}_2}\cap(\gS_{\ga'}\cap\gS_{\ti{\ga}})=
 \gC_{\ti{L}_2}\cap\gC_{L'_1\oplus\ti{L}_1}
\]
is also empty. Reasoning in the same way, for any subgroup
$H\subset\ti{L}_3$ of order $4$ such that $H\neq\ti{L_2}$ we have
$H$ is not a subset of ${L'}_1^\perp$. By Proposition \ref{Plev}
\[
 \#(D_{\ti{\ga}}^{-1}(H/\ti{L}_1)\cap
 D_{\ti{\ga}}^{-1}((L'_1\oplus\ti{L}_1)/\ti{L}_1))=1
\]
and therefore the cardinality of the intersection $\gC_H\cap\gS_{\ga'}$ is $2$.
\end{proof}

\begin{dsc}
We now move to the geometric part of the construction. Recall that
$C$ is a generic genus $3$ curve and that
$\mathcal{L}=(\langle \ga \rangle=L_1 \subsetneq L_2 \subsetneq
L_3)$ is a full isotropic flag in $\jac(C)[2]$. We start with the
Coble-Recillas construction. Next we construct the double cover
$Y\ra E$ mentioned above and review its properties. Finally we
describe the AGM construction through the bigonal construction.
\end{dsc}
\begin{ntt}
If $V \to U$ is a cover without specific name for the morphism, we
denote the morphism by $\pi_{V/U}$. In the following diagram we
summarize some constructions and notations which will be introduced
afterwords. Note that this diagram admit a symmetry with respect to
the vertical axis passing through the right term. These symmetric
objects, which are related to $(C',\mathcal{L}')$ under those
constructions, will be denoted with $'$. \label{Ddiagram}
\[
\xymatrix{& W \ar[dl] \ar[dr]^{/S} \\
 C \ar [ddrr] & & Z \ar [d]_{/\gs} \ar[rd]^{/i} & & M \ar[dl] \\
 & & X \ar [d]^{g^1_3} \ar[dr]^{/j} & Y \ar[d] & \\
& &|K_C+\ga|^* & E \ar[dr] &\\
& & & & \BP^1 }
\]
The left  construction is Coble-Recillas trigonal construction, and the right
construction is (half of) a bigonal construction (see \cite{Don} p.
68-69).
\end{ntt}
\begin{dsc}\label{Ctrigo}
 Let us recall some elements of the theory of the
trigonal construction (see for instance
\cite{Re},\cite{Dol1},\cite{Le}). Define
\[
  W:=\ol{C\times_{|K_C+\ga|^*}C\sm\gD_C}=\{(p_1,p_2)\in C\times C \; | \;
  p_1+p_2<K_C+\ga\}
\]
(where by $\gD_C$ we denote the diagonal). The curve $W$ admits a
natural involution - the coordinates switching which we denote by $S$. The
curve $Z=W/S$ can be viewed as the subset of $\Sym^2(C)$ defined
by $\{p_1+p_2 \; | \; p_1 +p_2 < K_C+\ga\}$. The curve $Z$ admits
three natural non-trivial involutions :
\begin{itemize}
\item $\gs : p_1+p_2 \mapsto p_3+p_4$ such that
$p_1+p_2+p_3+p_4 \sim K_C+\ga$; denote $X:=Z/\gs$.
\item $i : p_1+p_2 \mapsto p_3+p_4$ such that $p_3 + p_4
\sim p_1+p_2 +\ga$; denote $Y:=Z/i$.
\item $j=\gs \circ i$; denote $F:=Z/j$.
\end{itemize}
\begin{remark}\label{rem_theta}
The curves $Z$ has a ``theta divisor interpretation'': the Abel map
$\Sym^2C\ra\Pic^2C$ induces an isomorphism
$Z\cong\gTh_C\cap(\gTh_C+\ga)$ (see \cite{Le} 3.1-3.4). In
particular, $Z$ is generically a smooth curve of genus $7$. Note
that the involution $i$ is then $d \mapsto d +\ga$, $j$ is $d
\mapsto K_C-d$ and $\gs$ is $d \mapsto K_C+\ga-d$.
\end{remark}
Let us denote by $\psi$ the morphism  from $\Sym^2(C)$ to
$|K_C|$ defined by sending $p_1+p_2$ to the line $\ol{p_1p_2}$. As
the supports of $p_1+p_2$ and $j(p_1+p_2)$ as points on $C \subset
|K_C|^*$ lie on the same line, $\psi$ induces a morphism from
$F=Z/j$ to $|K_C|$ making the following diagram commutative:
\[
\xymatrix{ Z \ar@{^{(}->}[r] \ar[d]_{/j} & \Sym^2(C) \ar[d]^{\psi} \\
F \ar@{^{(}->}[r] & |K_C|. }
\]
It is classical that bitangents to $C$ are in $1-1$ correspondence
with odd theta characteristics. Moreover if $l$ is a bitangent
 corresponding to an element $\theta \in \gS_{\ga}$,  it
defines a point on $Z$ (still denoted $\theta$): indeed,
 if
  $l \cdot C=2(p_1+p_2)$ and if
  $\theta+\ga$ corresponds to a bitangent  with divisor $2(p_3+p_4)$
  one gets $(p_1+p_2)-(p_3+p_4) \sim \ga$, so $$p_1+p_2+p_3+p_4 \sim
  2(p_1+p_2)+\ga \sim K_C +\ga.$$
Thus $p_1+p_2 < K_C+\ga$. \\
By definition of $Z$ and $j$, the morphism $j$ is ramified exactly
at the $12$ elements of the Steiner system $\gS_{\ga}$. As $Z$ is
of genus $7$, $F$ is of genus $1$ and by the preceding embedding the
$12$ bitangents (viewed as points in $|K_C|$) are points on $F$. Note that
$F\subset|K_C|$ is a cubic: the degree of the map $F\ra|K_C|$ is $3$
since $\deg(\psi)=6$, and the map is non degenerate
because the points of $\psi(\gS_\ga)$ are not colinear.\\
As $i$ is fixed point free and commutes with $j$, it defines a fixed
point free involution $i_F$ on $F$ exchanging the points
$\theta,\theta+\ga \in \gS_{\ga}$. This involution defines a
point $\ga_F \in \Pic^0(F)[2]$ such that for all $p \in F$,
$p+\ga_F \sim i_F(p)$. The quotient of $F$ by the involution
$i_F$ is the curve $E=X/j$ (because $(Z/\gs)/j=(Z/j)/i$). Hence
$E$ is naturally embedded  in $|K_C|^*$ by the image of $\pi_{F/E} :
p \mapsto p\cap (p+\ga_F)$.
\end{dsc}
\begin{theorem}[Coble, see \cite{Co} sections 47-49,\cite{Le}]\label{TPC}
The images of the points of $\gS_\ga\subset F$ under the map
$\pi_{F/E}$ sit on a unique conic $Q\subset|K_C|^*$. They are the
intersection points of the bitangents $\theta,\theta+\ga \in
\gS_{\ga}$. The locus $Q \cap E$ is the ramification locus of the
map $\pi_{Y/E} : Y \to E$.
\end{theorem}
\begin{dsc}[The folloing results are mostly due to Coble and to Recillas]\label{Dinverse}
Let  $\ga_E$ be  the unique non zero element in
$(\pi_{F/E})_{*}(\Pic(F)[2])\subset\Pic(E)[2]$. Given the double
cover $Y\ra E$ and the $2$-torsion point $\ga_E$ one can reconstruct
the curve $C$ and the linear system $|K_C+\ga|$ in the following
way. We have $F \simeq E/\ga_E$. The curve $Z$ is isomorphic to
the fibered product $Y\times_E F$. This construction induces two
commuting involutions $i,j$ on $Z$. Since the double cover
$\pi_{F/E}$ is unramified, the involutions $i,j$ on $Z$ are fixed
points free. The genera of the curves $Z,X,Y$ are then $7,4,4$
respectively. The curve $X$ is a bielliptic curve of genus $4$ which
has only one $g^1_3$ up to the bielliptic involution.  Thus we are back in the
trigonal construction setting.
Note that $C$ and the linear system $|K_C+\ga|$ are invariant under the
choice of the $g^1_3$.
On our way to $(C',\mathcal{L}')$, we have now expressed
$(\jac(C),\ga)$ in terms of the double cover $Z/X$. The second step in the
construction (see Theorem \ref{TCYE}) is to interpret the symplectic data
through the quotient $i$. Our main tool for analysis of non
principally polarized Abelian varieties are Lemmas \ref{LDLisogeny}
and \ref{mono} below.
\end{dsc}
\begin{lemma}[see \cite{DL}, Lemma 1]\label{LDLisogeny}
Let $\ti{V}\ra V$ be an admissible double cover, and let
 $\nu\ti{V}\ra \nu V$ be
its partial normalization at $r>1$ points $x_1,\ldots,x_r\in V$. Let g be the
arithmetic genus of the partial normalization $\nu V$, so the arithmetic
genus of $V$ is $g+r$. Then $\prym(\ti{V}/V)$ has a principal polarization,
$\prym(\nu\ti{V}/\nu V)$ has a polarization of type $2^g 1^{r-1}$, and the
pullback map
\[
  \nu^* : \prym(\ti{V}/V)\ra\prym(\nu\ti{V}/\nu V)
\]
is an isogeny of degree $2^{r-1}$.
\end{lemma}
\begin{lemma}[the monodromy argument] \label{mono}
Let $V'\ra V$ be a finite cover such that the
Galois group of the Galois closure of $V'/V$ is $2$-transitive on the cover,
then the only section of the cover $V'\times_V V'\ra V$ is the diagonal.
\end{lemma}
\begin{dsc}
We will apply lemma \ref{mono} with the covers $\mathcal{A}_3^1
\to \mathcal{A}_3$ and $\mathcal{A}_3^F \to \mathcal{A}_3$ where
$\mathcal{A}_3^1$ (resp. $\mathcal{A}_3^F$) is the moduli space of
PPAV's dimension $3$  with a $2$-torsion point
(resp. with a maximal isotropic group).
\end{dsc}
\begin{theorem}\label{TCYE}
The quotient by $i$ induces an isogeny of Abelian varieties:
\[ \phi :
  \prym(Z/X)\ra \prym(Y/E).
\]
Identifying the principally polarized Abelian varieties $\jac(C)$
and $\prym(Z/X)$ as in \cite{Don} Theorem 2.11 (p. 76), the kernel
of $\phi$ is $\ga^\perp$.
\end{theorem}
\begin{proof}
The proof consists of three steps:
\item
{\em Step 1: The map $\phi : \prym(Z/X)\ra \prym(Y/E)$ is an isogeny :}\\
To prove this claim it suffices to prove that the induced map on the tangent
spaces at $0$ is an isomorphism.  We do this by considering the
space $M:=H^0(Z,\gO^1_Z)$ as a $\gal(Z/E)$ module,
and calculating the module decomposition to irreducible representations.
We denote by $M_{-}$ the irreducible representation corresponding to the
character whose kernel is the subgroup $\langle -\rangle\subset \gal(Z/E)$,
and by $M_1$ the irreducible representation corresponding to the trivial
character. Using these notations we have:
\[
  \begin{aligned}
    H^0(Z,\gO^1_Z)&=M_1 \oplus M_i \oplus M_j\oplus M_\gs,
     \quad H^0(E,\gO^1_E)=M_1 \\
    H^0(F,\gO^1_F)&=M_1^j\oplus M_i^j\oplus M_j^j\oplus M_\gs^j=
      M_1\oplus M_j,\\
    H^0(Y,\gO^1_Y)&=M_1^i \oplus M_i^i \oplus M_j^i \oplus M_\gs^i=
     M_1\oplus M_i, \\
    H^0(X,\gO^1_X)&=M_1^{\gs}\oplus M_i^{\gs}\oplus M_j^{\gs}\oplus M_\gs^{\gs}
      =M_1\oplus M_{\gs}.
  \end{aligned}
\]
However, since $E$ and $F$ are both of genus $1$,
the map $H^0(F,\gO^1_F)\ra H^0(E,\gO^1_E)$ is an isomorphism.
Thus we have $M_j=0$, and our claim holds.
\item {\em Step 2: The kernel of the isogeny $\phi$ is a subset of $\prym(Z/X)[2]$:}\\
Denoting by $[2]_-$ the multiplication by $2$ on an Abelian variety we have
\[
  {\pi_{Z/Y}}_*\pi_{Z/Y}^*:\jac(Y)\ra\jac(Y)=[2]_{\jac(Y)}.
\]
Denote by $\mu_*,\mu^*$ the restrictions of maps
${\pi_{Z/Y}}_*,\pi_{Z/Y}^*$ to the Abelian varieties
$\prym(Y/E),\prym(Z/X)$ respectively. Note that $\phi=\mu^*$. We
have $\mu^*\mu_*=[2]_{\prym(Z/X)}$.
\item {\em Step 3: Computation of the kernel of $\phi$:}\\
By applying Lemma \ref{LDLisogeny} to some degeneration of the cover
$Y/E$ along its ramification locus $Q \cap E$, on finds that the
polarization type of the variety $\prym(Y/E)$ is $2^1 1^5$. Thus,
the order of
the kernel of $\phi$ is $32$.\\
By \cite{Don} Theorem 2.11 (p. 76) the norms in the trigonal
construction induce an isomorphism $\jac(C)\cong\prym(Z/X)$. The
kernel of $\phi$ is can thus be identified with $\gb^\perp$ for some
$\gb\in\jac(C)[2]$. The map $(\jac(C),\ga) \to (\jac(C),\gb)$ gives a
endomorphism of $\mathcal{A}_3^1$ and then a section from
$\mathcal{A}_3^1$ to $\mathcal{A}_3^1 \times_{\mathcal{A}_3}
\mathcal{A}_3^1$. By Lemma \ref{mono} this section maps into the
diagonal component so $\ga=\gb$.
\end{proof}
\begin{dsc}\label{Dlev}
Denote by $q_1,\ldots,q_6$ the intersection points of $E$ and $Q$.
In Proposition \ref{Plev} we identified pairs of $q_i$s with the
non-zero points of the symplectic space $\ga^\perp/\ga$. This
identification induces bijections between the following sets of
data:

\ \\
\begin{tabular}{l|l|l}
data on isotropic subgroups in & data on isotropic  & partitions of the \\
$\jac(C)[2]$ that contain $\ga$. & subgroups of $\ga^\perp/\ga$. &  points $\{q_i\}_{i=1\ldots 6}$\\
\hline
\hline
Iso. subgps of order $4$ & Iso. subgps of order $2$ & $2+4$ \\
\hline
Maximal iso. subgps & Maximal iso. subgps. &  $2+2+2$ \\
\hline
Full isotropic flags & Full isotropic flags & $2+(2+2)$\\
\end{tabular}
\ \\ \\

We let $\{q_1,q_2\},\{\{q_3,q_4\},\{q_5,q_6\}\}$ be the partition of
the $q_i$s that corresponds to the full isotropic flag $\caL$.
Denote by $\pi_{E/\BP^1}$ the linear system $|q_1+q_2|$ on the curve
$E$ and by $B$ the ramification locus of $\pi_{E/\BP^1}$. Note that
the symmetric construction introduces a set $B'$. Using these
definitions we are ready to prove the correctness of our
construction:
\end{dsc}
\begin{theorem}\label{TAGM}
Let assume that the ramification pattern of the tower $Y\ra E \ra
\BP^1$ is generic. Denote by $\ti{H}\ra H \ra \BP^1$ the image of
the tower $Y\ra E\ra \BP^1$ by the bigonal construction. Then the
tower $Y'\ra E' \ra \BP^1$ is the normalization of the tower
$\ti{H}\ra H \ra \BP^1$. Moreover, there is a 1-1 correspondence
between points $b\in Q\cap E\sm\{q_1,q_2\}$ and points $b'\in B'$,
given by
\[
  \pi_{E/\BP^1}(b)=\pi_{E'/\BP^1}(b').
\]
\end{theorem}
\begin{proof}
 The ramification pattern
of the bigonal construction on $Y\ra E\ra\BP^1$ is the following
(see \cite{Don} p. 68-69):
\begin{itemize}
\item If $\pi_{E/\BP^1}^{-1}(a)=\{q_1,q_2\}$, then $Y/E$
is ramified over both $q_1,q_2$, $\pi_{H/\BP^1}^{-1}(a)$ is a node,
and over this node the curve $\ti{H}$ is a gluing of two ramified
sheets (symbolically : $\subset \subset / = | \supset \! \subset /
\times$) .
\item If $\pi_{E/\BP^1}^{-1}(a)=2b$ for some $b\in B$ then $Y/E$ is \'etale over $b$,
$\pi_{H/\BP^1}$ is \'etale over $a$, and $\ti{H}$ is ramified over
one of the points in $\pi_{H/\BP^1}^{-1}(a)$ and ramified over the
other one (symbolically : $\subset \subset / \subset | \subset = /
=$).
\item If $\pi_{E/\BP^1}^{-1}(a)\ni q$ for some $q\in Q\cap E\sm\{q_1,q_2\}$ then
$Y/E$ is ramified over $q$, and \'etale over the other point in
$\pi_{E/\BP^1}^{-1}(a)$. Moreover, $\pi_{H/\BP^1}$ is ramified at
$a$, and $\ti{H}/H$ is \'etale over both branches (symbolically :
$\subset = / = | \subset \subset / \subset$).
\item In all other points the ramification patterns of the towers
$Y/E/\BP^1$ and $\ti{H}/H/\BP^1$ are generic (i.e. unramified).
\end{itemize}
 Denote by $\nu\ti{H},\nu H$ the normalizations of the curves
$\ti{H},H$ respectively. By the Riemann-Horowitz formula, the genera
of the curves $\nu\ti{H},\nu H$ are $4,1$ respectively. The
ramification pattern over the points
$\pi_{E/\BP^1}(q_3),\pi_{E/\BP^1}(q_4)$,
$\pi_{E/\BP^1}(q_5),\pi_{E/\BP^1}(q_6)$ in the tower $Y\ra
E\ra\BP^1$ is $\subset =/=$. Thus, the partition
$\{\{q_3,q_4\},\{q_5,q_6\}\}$ induces a partition to two pairs of
the $4$ ramification points of the map $\nu H\ra\BP^1$, which
induces a choice of a $2$-torsion point in $\Pic^0(\nu H)$. Applying
the reconstruction technique from \ref{Dinverse} to the double cover
$\nu\ti{H}\ra \nu H$ and the $2$-torsion point we get a smooth curve
$C''$ of genus $3$. We claim that we have the following degrees for
the isogenies
$$
\jac(C)  \overset{2}{\to} \prym(\ti{H}/H) \overset{2}{\to}
\prym(\nu\ti{H}/\nu H) \overset{2}{\to} \jac(C'').$$  By \cite{Pa}
Proposition 3.1 (page 307) the Abelian variety $\prym(\ti{H}/H)$ is
isomorphic to the dual of the Abelian variety $\prym(Y/E)$. Since
$\prym(Y/E)$ is isomorphic to $\jac(C)/\ga^\perp$, we have by
\ref{DPPAV} that $\prym(\ti{H}/H) \simeq \jac(C)/\ga$. The second
arrow is a consequence of Lemma \ref{LDLisogeny} for the
normalization of $\ti{H}/H$ over the point of type $ \supset \!
\subset / \times$. The third arrow
follow from Theorem \ref{TCYE}.\\
Thus we have obtained an isogeny of degree $2^3$. By our observation
in \ref{DPPAV}, this isogeny is given by a maximal isotropic group
$L \in \jac(C)[2]$. In the same spirit as in the proof of \ref{TCYE}
we can use the monodromy argument of Lemma \ref{mono} for $L_3$ and
$L$ to prove that $L=L_3$. Thus $C'' \simeq C'$ and we get
$$
\nu\ti{H}\simeq Y',\quad \nu H\simeq E'.
$$
The ramification patterns prove the last assertion of the theorem.
\end{proof}
\begin{remark}
In the same way, one can describe analogous results for the other
non-generic ramification patterns of the tower $Y \ra E\ra\BP^1$,
but this effort is redundant: in Section \ref{Sprog} we will find a
formula which gives $C'$ in terms of a generic pair $(C,L)$; since
$C'$ is continuous in the pair $(C,L)$, the formula will be correct
for all pairs $(C,L)$ for which the denominators in the formula are
non zero.
\end{remark}
%
\section{The isomorphism between the canonical classes}\label{Scanon}
%
\begin{dsc}\label{Dcanon}
 In this section we describe the
isomorphism $k:|K_{C'}|^*\ra|K_C|^*$ between the duals of the
canonical classes of the curves $C$ and $C'$. In Section
\ref{Sprog}, this description is used to calculate the equation of
the canonical embedding of curve $C'$ in terms of the canonical
embedding of the curve $C$ and to describe $\ol{k}: H^0(K_C) \to
H^0(K_{C'})$.

We describe the isomorphism $k$ by considering the images and
preimages (under the map $k$) of the sets $B,B'$ (recall the
definition in \ref{Dlev}) and the points defined below:
\[
  p:=E\cap\ol{q_1q_2}\sm\{q_1,q_2\},\quad
  p':=E'\cap\ol{q'_1q'_2}\sm\{q'_1,q'_2\}.
\]
The resulting description is encoded in the following theorem:
\end{dsc}
\begin{theorem}\label{TagM}
The isomorphism $k$ is completely determined by the following identities:
\[
  \begin{aligned}
   k(Q'\cap E'\sm\{q'_1,q'_2\})=B, \quad&k(\ol{q'_1q'_2})=\ol{q_1q_2}\\
   k(B')=Q\cap E\sm\{q_1,q_2\},\quad&
  k(p')=p,
   \end{aligned}
\]
where the identification of $Q'\cap E'\sm\{q'_1,q'_2\}$ with $B$, and of
 $B'$ with $Q\cap E\sm\{q_1,q_2\}$ are the ones from Theorem \ref{TAGM}.
\end{theorem}
\begin{proof}\label{PFTagM}
The theorem follows from Theorem \ref{Tpsi} and two applications of
Theorem \ref{Tphi} below.
\end{proof}
\begin{dsc}\label{Dkdesc}
To describe the isomorphism $k$ we present it as a composition of
three isomorphisms. Denote by $j_Y$ (resp. $J_Y$) the involution on
the curve $Y$ (resp. the homology group $H^0(K_Y)$) induced from the
double cover $Y\ra E$. Denote by $H^0(K_Y)_{\odd}$ (resp.
$H^0(K_Y)_{\even}$) the odd (resp. even) part of $H^0(K_Y)$ with
respect to the involution $J_Y$. We denote by $|K_Y|_{\odd}$
(respectively $|K_Y|_{\even}$) the projectivization of the vector
space $H^0(K_Y)_{\odd}$ (respectively $H^0(K_Y)_{\even}$). We use
the analog notations for subspaces of $|K_{Y'}|$. The involution
$J_Y$ induces an involution on the dual of the canonical system
$|K_Y|^*$. The fixed set under this involution is the union of the
projective plane $|K_Y|_\odd^*$ and a point $p_Y$,
the projectivization of the space $H^0(K_Y)_\even$.\\
 By Theorem \ref{TAGM} we have a
sequence of isogenies of Abelian varieties
\[
  \jac(C)\overset{/\ga^\perp}{\ra}\prym(Y/E)\ra\prym(Y'/E')
  \overset{/{\ga'}^\perp}\longleftarrow\jac(C').
\]
Taking the tangents spaces at $0$ of these varieties we get a sequence of
isomorphisms:
\begin{equation}\label{Etangents}
  H^0(K_C)\overset{\ol{\phi}}{\ra}H^0(K_Y)_{\odd}\overset{\ol{\psi}}{\ra}H^0(K_{Y'})_{\odd}
  \overset{\ol{\phi'}}\longleftarrow H^0(K_{C'}),
\end{equation}
 Taking the dual, inverse and projectivizations of the spaces and morphisms
in Equation (\ref{Etangents}) we get another sequence of
isomorphisms:
\[
  |K_C|^*\overset{\phi}{\longrightarrow}|K_Y|_{\odd}^*\overset{\psi}{\longrightarrow}
  |K_{Y'}|_{\odd}^*
  \overset{\phi'}{\longleftarrow} |K_{C'}|^*.
\]
By construction, these two morphisms have interpretation in terms of
the trigonal and bigonal constructions. With the notations of
\ref{Ctrigo} the morphism $\phi$ is induced by ${\pi_{W/Y}}_*
\pi_{W/C}^*$. In the same way, denoting by $M$ the normalization of
the Galois closure of the tower $Y\ra E\ra\BP^1$, the isomorphism
$\psi$ is defined as the composition
${\pi_{M/Y'}}_*\circ\pi_{M/Y}^*$.\\
We discuss the morphism $\phi$ in Theorem \ref{Tphi} below. Our
analysis is based on the two views of the set $\{q_i\}_{i=1 \ldots
6}$ presented in Theorem \ref{TPC}:
\begin{itemize}
\item The $q_i$s are in natural 1-1 correspondence with intersection points
of pairs of bitangents, which lie in $|K_C|^*$.
\item The $q_i$s are in natural 1-1 correspondence with the fixed points in
$Y$ of the involution $j_Y$, i.e. with the points
of $Y\cap|K_Y|_\odd^*\subset|K_Y|^*$.
\end{itemize}
We interpret the relation between these views using the norms in the
trigonal construction. If $\pi_{V/U} : V\ra U$ is one of the covers
arising in our construction, we denote  the ramification divisor of
the map $\pi_{V/U}$
 by $R_{V/U}$.
We will make repetitive use of the following version of the Riemann-Horowitz
theorem (see \cite{Har} Proposition IV.2.1):
{\em Let $\gom$ be a differential on $U$ then the zero divisor of the differential $\pi_{V/U}^*\gom$ is $\pi_{V/U}^*((\gom)_0)+R_{V/U}$}.
\end{dsc}
\begin{theorem}\label{Tphi}
The map $\phi :|K_C|^* \ra|K_Y|^*_\odd$ takes each of the $q_i$s to
the corresponding point in $Y\cap|K_Y|_\odd^*\subset|K_Y|^*$.
Moreover, this property defines $\phi$.
\end{theorem}
\begin{proof}
To avoid confusion between the points of $W,Z,Y$ and divisors, we
denote here a point of $W$ by $(p_1,p_2)$, on $Z$ by  $(p_1+p_2)$
and
a point of $Y$ by $\{(p_1+p_2),(p_3+p_4)\}$ if $p_3+p_4=j(p_1+p_2) \in Z$.\\
Let $\{p_1+p_2,p_3+p_4\}$ be the pair of theta characteristics in
$\gS_\ga$ which represent one of the $q_i$s (see \ref{Dlev}).
Let $\gom\in H^0(K_C)$ be a differential such that
$\gom_0=2(p_1+p_2)$ then
\[
  (\pi_{W/C}^*(\gom))_0-R_{W/C}=\pi_{W/C}^*(\gom_0)=2((p_1,p_2)+(p_2,p_1)
             +\sum_{\substack{1 \leq i \leq 2\\ 3 \leq j \leq 4}}(p_i,p_j)).
\]
So $(\pi_{W/C}^*(\gom)) \geq 2((p_1,p_2)+(p_2,p_1)) + R_{W/C}$.
 By the properties of the trigonal construction (see \cite{Don} p. 74)
we have $R_{W/C} \geq R_{W/Z}$. So we get
$$
(\pi_{W/C}^*(\omega)+S\pi_{W/C}^*(\gom))_0\geq2((p_1,p_2)+(p_2,p_1))+R_{W/Z}.$$
The left summand is precisely the pull back of
${\pi_{W/Z}}_*\pi_{W/C}^*(\gom)$ so finally get
$$
 ({\pi_{W/Z}}_*\pi_{W/C}^*(\gom))_0\geq2 (p_1+p_2).
$$
Now
\[
{\pi_{Z/Y}}_*{\pi_{W/Z}}_*\pi_{W/C}^*(\gom)=
{\pi_{W/Y}}_*\pi_{W/C}^*(\gom)= \ol{\phi}(\gom)
  \in H^0(K_Y)_\odd.
\]
So by invariance under $J_Y$
$$(\ol{\phi}(\omega))_0=
  ({\pi_{W/Y}}_*\pi_{W/C}^*(\gom))_0  \geq
      2 \{(p_1+p_2),(p_3+p_4)\}.
$$
As $q_i$ is the intersection points of the two bitangents supported by
$(p_1,p_2),(p_3,p_4)$, the last inequality implies
\[
  \phi(q_i)=\{(p_1+p_2),(p_3+p_4)\}.
\]
Since the $q_i$s are $6$ non-collinear points this property
completely describes the map $\phi$.
\end{proof}
\begin{dsc}
It remains to analyze the isogeny $\psi$ induced from the bigonal
construction relating the double covers $Y\ra E$ and $Y'\ra E'$ (see
the proof of Theorem \ref{TAGM}). We make the identifications of
Theorem \ref{Tphi} (so $\psi$ becomes $k$). As the linear system
$|p+q_1+q_2|$ spans the space $|K_Y|_{\odd}^*$ (and the same with
the symmetric notation $'$) we reduce the description of the map
$\BP\psi$ to a description of a natural isomorphism between these
linear systems on the curves $E,E'$.
\end{dsc}
\begin{theorem}\label{Tpsi}
The isomorphism $\psi$ is determined by the following identities:
\[
  \begin{aligned}
   \psi(Q'\cap E'\sm\{q'_1,q'_2\})=B, \quad& \psi(\ol{q'_1q'_2})=\ol{q_1q_2}\\
   \psi(B')=Q\cap E\sm\{q_1,q_2\},\quad&
  \psi(p')=p,
   \end{aligned}
\]
where the identification of $Q'\cap E'\sm\{q'_1,q'_2\}$ with $B$, and of
 $B'$ with $Q\cap E\sm\{q_1,q_2\}$ are the ones from theorem \ref{TAGM}.
\end{theorem}
\begin{proof}
The theorem follows from the two claims below:
\begin{enumerate}
\item Let $t$ be a point in $\BP^1$ and let $\gom$ be a differential in
$H^0(K_Y)_\odd$ such that
$\half{\pi_{Y/E}}_*((\gom)_0)=p+\pi_{E/\BP^1}^*(t)$ then the image
of $\gom$ under $\ol{\psi}$ satisfies
$\half{\pi_{Y'/E'}}_*((\ol{\psi}(\gom))_0)=p'+\pi_{E'/\BP^1}^*(t)$.
\item Let $b$ be a point in the set $B\subset E$ and let $q'_i$ be the
corresponding point (in the sense of Theorem \ref{TAGM}) in the set
$Q'\cap E'$. Let $\gom$ be a differential in $H^0(K_Y)_\odd$ such
that ${\pi_{Y'/E'}}_*(\gom)_0\geq 2b$ then the image of $\gom$ under
$\ol{\psi}$ satisfies $(\ol{\psi}(\gom))_0\geq 2q'_i$.
\end{enumerate}
As in the proof of Theorem \ref{Tphi} we make repetitive use of
Riemann-Horowitz theorem.

\item {\em Proof of claim 1:}
Let $\gom$ be a differential as in the first claim above.
Since the zero divisor of the differential $\gom$ is
moving with $t$, the intersection
$(\pi_{M/Y}^*(\gom))_0\cap R_{M/Y}$ is generically empty.
By the definition of the bigonal construction
\[
  \begin{aligned}
  (\ol{\psi}(\gom))_0=&({\pi_{M/Y'}}_*\pi_{M/Y}^*(\gom))_0=
  \half{\pi_{M/Y'}}_*((\pi_{M/Y}^*(\gom))_0-R_{M/Y})\\
  =& \half{\pi_{M/Y'}}_*\pi_{M/Y}^*({\pi_{Y/\BP^1}}^*(t)+{\pi_{Y/E}}^*(p)).
  \end{aligned}
\]
Then
$$\frac{1}{2} {\pi_{Y'/E'}}_*((\ol{\psi}(\omega))_0)=\frac{1}{2} ({\pi_{M/E'}}_*
\pi^*_{M/\BP^1}(t)+ {\pi_{M/E'}}_* \pi^*_{M/E}(p))=\pi^*_{E/\BP^1}(t)+p'.$$
\item {\em Proof of claim 2:} Define the objects $\gom,b,q'_i$ as in
the second claim above. By the definition of the point $b$ we have
$(\pi_{M/Y}^*(\gom))_0\geq \pi_{M/\BP^1}^*(\pi_{E/\BP^1}(b))$. Whence:
\[
  \begin{aligned}
  (\ol{\psi}(\gom))_0=&({\pi_{M/Y'}}_*\pi_{M/Y}^*(\gom))_0=
  \half{\pi_{M/Y'}}_*((\pi_{M/Y}^*(\gom))_0-R_{M/Y'})   \\
   \geq&\half{\pi_{M/Y'}}_*(\{t\in\pi_{M/\BP^1}^*(\pi_{E/\BP^1}(b))|t\not\in R_{M/Y'}\}).
  \end{aligned}
\]
To prove the inequality $(\ol{\psi}(\gom))_0\geq 2q_i$ it suffices
to show that $R_{M/Y'}\cap\pi_{M/Y'}^*R_{Y'/E'}=\emptyset$. By the
bigonal construction dictionary (see Theorem \ref{TAGM}), if the
cover $Y\ra\BP^1$ is ramified over a point $t$ then the cover
$E'\ra\BP^1$ is \'etale over $t$. Since the curve $M$ can be defined
as the product $E'\times_{\BP^1}Y$ there are no multiple points in
the ramification divisor $R_{M/\BP^1}$. This proves that
$R_{M/Y'}\cap\pi_{M/Y'}^*R_{Y'/E'}=\emptyset$.
\end{proof}
%
\section{A small matter of programming}\label{Sprog}
%
\begin{dsc}\label{Dprog}
In the previous sections we presented an explicit construction of
the AGM in genus $3$. In this section we close the gap between
``explicit'' and a formula. We tackle five problems: describing the
pair $(C,\ga)$ (\ref{Data1}-\ref{eqC}), describing the
intermediate data $(E,Q),(E',Q')$ in term of $(C,L_2)$
(\ref{Tcoor}-\ref{Data14}), describing the pair $(C',\caL')$ in
terms of the pair $(C,\caL)$ (\ref{DcalcAis}-\ref{TcalcAis}),
describing the isomorphism $\ol{k}:H^0(C,K_C)\ra H^0(C',K_{C'})$
(\ref{Cdiff}), and describing a flag $\ti{\caL}$ (see \ref{Dbasis})
in terms of $E',Q',\caL'$ (\ref{Dcalcbasis}-\ref{Dexpbasis}). We
solve these problems by a ``coordenification'' of the proof of
Theorem \ref{TagM}. We identify the two spaces $|K_C|^*$ and
$|K_{C'}|^*$ under the isomorphism $k$. We denote the coordinates on
this space by $x,y,z$ (a precise choice of coordinate is described
in \ref{Data1}).
\end{dsc}
\begin{ntt}
To write the equations, we use the lexicographic order on the dual
coordinates of $x,y,z$. i.e. instead of writing
$ax^2+bxy+cxz+dy^2+eyz+fz^2$, we will write $(a,b,c,d,e,f)$. If we
will talk about the corresponding curve, we will use the projective
coordinates $(a;b;c;d;e;f)$. Finally, we will abuse the notations
by using the name of a curve in $\BP^2$ for its defining equation.
\end{ntt}

\begin{dsc}{\em The data $(C,\ga)$.} \label{Data1}
 Let us consider the natural bilinear map
\[
  \begin{aligned}
  H^0(K_C+\ga)\times H^0(K_C+\ga)&\ra H^0(2K_C)=H^0(\cO_{|K_C|}(2))\\
((a s_1 + b s_2),(c s_1+ d s_2))&\mapsto acA_1+(bc+ad)A_3+bd A_2
  \end{aligned}
\]
defined by the tensor multiplication of the sections $s_i$.
Identifying $|K_C+\ga|$ with $\BP^1$ we also get a map
\[
  m:\BP^1 \times \BP^1 \to |2 K_C|.
\]
We simplify the conics $A_i$ by a special choice of our
sections $s_i$: Let $\theta_i,\theta_i+\ga \in \gS_{\ga}$ for
$i=1,2$; we assume below that the pair
$(\theta_i,\theta_i+\ga)$ corresponds to the points $q_i$ by
Theorem \ref{TPC}. Note that the choice of a distinguished pair
$\{q_1,q_2\}$ is
equivalent to the choice of $L_2$ in $\jac(C)[2]$ as described in \ref{Ctrigo}.\\
Let  $l_{1i}$ (resp. $l_{2i}$)  the bitangents corresponding to
$\gth_i$ (resp. $\gth_i+\ga$) and let be $D_{ji}$ the effective
divisor of degree $2$ such that $2 D_{ji}=(l_{ji})_0$, the we have
\[
  K_C+\ga \sim D_{1i}+\gs(D_{1i}) \sim D_{1i}+i\circ j(D_{1i}) \sim
D_{1i} + i(D_{1i}) \sim D_{1i} + D_{2i}.
\]
We denote by $s_i$ the
sections of $H^0(K_C+\ga)$ corresponding to $D_{1i}+D_{2i}$. With this
choice of sections, we have $A_1=l_{11} l_{21}$
and $A_2=l_{12} l_{22}$. We now fix the coordinates $(x;y;z)$ of
$|K_C|^*$ such that $A_1=(y-z)(y+z)$ and $A_2=(x-z)(x+z)$.
The following Proposition is now a particular case of a classical result (see
\cite{Dol2}).
\begin{proposition} \label{eqC}
The quartic $C$ is given by $A_3^2-A_1 A_2=0$. Assume that the
quadrics $A_1,A_2,A_3$ are respectively given by
\[
  (0,0,0,1,0,-1),\quad (1,0,0,0,0,-1),\quad
  (a,b,c,d,e,f),
\]
\\
then the coordinates of $C \in |\mathcal{O}_{|K_C|^*}(4)|$ are
\[
  \begin{aligned}
  (&a^2; 2ab; 2ac; b^2 + 2ad-1; 2bc + 2ae;1+c^2+2af;2bd;2cd+2be;2ce+2bf;\\
  ;&2cf;d^2; 2de; 1 + e^2 + 2df; 2ef;f^2-1).
  \end{aligned}
\]
\end{proposition}
\noindent In \ref{Tcoor}-\ref{Data14}, we will find the equations
for the intermediate data $(E,Q),(E',Q')$.
\end{dsc}
\begin{theorem}\label{Tcoor}
Assume that the quadrics $A_1,A_2,A_3$ are as in \ref{eqC}. The
coordinates of the curves $E,E'\in|\cO_{|K_C|^*}(3)|$ are
\[
  \begin{aligned}
  E=&(0;c;b;e;2(a+d+f);e;0;b;c;0),\\
  E'=&(-2ac;-2ae;b^2-c^2-4af-1-4a^2;2cd;4b(d-a);2be-4ac-2cf\\
     ;&2de;1-b^2+4d^2+e^2+4df;2e(2d+f)-2bc;e^2-c^2).
  \end{aligned}
\]
The coordinates of the conics $Q,Q' \in|\cO_{|K_C|^*}(2)|$ are
\[
  \begin{aligned}
  Q=&(0;ce(b^2-1+4ad)-2b(c^2d+e^2a)\\
     ;&b(-2bcd-e+b^2e)+2a(b^2c+ce^2-2bef)-4a^2be;0\\
     ;&b^3c+2c^2de+2b^2(d-a)e-bc(1+4d^2+4df)\\
     ;&c^2e^2+b^2(c^2+e^2)-2bce(a+d+f)),\\
  Q'=&(-a(e^2-c^2);0;(d-a)(c(a+d+f)-be);d(c^2-e^2)\\
     ;&2(a-d)(e(a+d+f)-bc);(d-a)(c^2-e^2)).
  \end{aligned}
\]
\end{theorem}

\begin{dsc}\label{Dcoor}
The essence of the proof is to convert the problem to a sequence of
``steps'' of the following form: {\em find a pencil that is spanned
by two known forms, and that contains another form we have to
calculate}. We perform these steps explicitly using a computer. Let
us stress though, that in most cases two of the three forms involved
in the computation are divisible by a known linear form. Thus, the
obstinate reader could still check the computations below, up to and
including Theorem \ref{TcalcAis}, by hand.\\
\item {\em The curve $E$.} Recall (see \ref{Ctrigo}) that to any
point $p_1+p_2 \in Z$ one associate a point $q \in E \subset
|K_C|^*$ as $\ol{p_1p_2} \cap \ol{p_3p_4}$ where $i(p_1+p_2)=
p_3+p_4 \in Z$. Denote by $A$ the conic given by the product of the
lines $\ol{p_1p_2}$ and $\ol{p_3p_4}$. Note that the singular point
of $A$ is $q$.
\begin{lemma} \label{Enet}
There exists a unique couple of sections (up to permutation) $(t_1,t_2) \in
H^0(K_C+\ga)^2$ such that $t_1 t_2=A$. Conversely, the singular
point of each singular conic in  $m(\BP^1 \times \BP^1)$ is on $E$.
\end{lemma}
\begin{proof}
Let us show the first assertion. By definition of $j$, the zero
divisor of $A$ is $2 K_C \sim p_1+p_2+ j(p_1+p_2) + p_3+p_4+
j(p_3+p_4).$ As $p_3+p_4=i(p_1+p_2)$, one gets
$$(A)_0=p_1+p_2+ j\circ i (p_1+p_2) + i(p_1+p_2)+ \gs \circ
i(p_1+p_2).$$ As $p_1+p_2+ \gs(p_1+p_2) \in |K_C+\ga|$ (resp.
$i(p_1+p_2)+ \gs \circ i(p_1+p_2) \in |K_C+\ga|$), the divisor
defines a unique section $t_1$ (resp. $t_2$) in $H^0(K_C+\ga)$.
\end{proof}
To account for the permutation of the two lines in Lemma \ref{Enet}, we
introduce the  map
\begin{eqnarray*}
  v_2 : & \BP^1 \times \BP^1 &\to  \BP^2 \simeq \Sym^2\BP^1\\
        & (\gl,\mu),(\gl',\mu') & \mapsto (\gl \gl',\mu\mu',\gl\mu'+\mu\gl').
  \end{eqnarray*}
Let $(X;Y;Z)$ be
the coordinates in $v_2(\BP^1 \times \BP^1)$. Lemma \ref{Enet} can
be reformulated into : the curve $E$ can be seen  as the locus
$(x;y;z)$ of singular points in the net of conics $X A_1+ Y A_2 + Z
A_3$. If $M=X_0 A_1+ Y_0 A_2 + Z_0 A_3$ is such a conic, it is
singular at $q_0$ if and only if
\[
\left(\begin{array}{c} 0 \\ 0 \\ 0 \end{array} \right)
=\left(\begin{array}{c} M_x(q_0) \\ M_y(q_0) \\M_z(q_0)\end{array}\right)
=\left(\begin{array}{ccc}
(A_1)_x(q_0) & (A_1)_y(q_0) & (A_1)_z(q_0)\\
(A_2)_x(q_0) & (A_2)_y(q_0) & (A_2)_z(q_0)\\
(A_3)_x(q_0) & (A_3)_y(q_0) & (A_3)_z(q_0)\\
\end{array}\right)^t \cdot
\left( \begin{array}{c} X_0 \\ Y_0 \\ Z_0 \end{array} \right).
\]
Denote by $\jac(A_1,A_2,A_3)$ the previous matrix. Thus the curve
$E$ is given by $\det(\jac(A_1,A_2,A_3))=0$.\\

\item {\em The curve $Q$.}  It is easy to check algebraically that the point
$o:=(0;0;1)$ lies on $E$ (see \ref{pointo} for a geometric explanation).
Let $\hat{q_1}:=\ol{oq_1}\cap E\sm\{o,q_1\}$ and
$\hat{q_2}:=\ol{oq_2}\cap E\sm\{o,q_2\}$. Let
\[
  \hat{Q}:=
 \mathrm{Nulls}(4\frac{\partial A_3}{\partial x}\frac{\partial A_3}{\partial y}
-\frac{\partial A_1}{\partial y}\frac{\partial A_2}{\partial x}).
\]
\end{dsc}
\begin{proposition}\label{PcalcQY}
The cubics $E,Q\ol{\hat{q_1}\hat{q_2}},\hat{Q}\ol{q_1q_2}$ lie on
the same pencil. Moreover $p \in \ol{\hat{q_1}\hat{q_2}}$.
\end{proposition}
\hspace{3cm}
\begin{picture}(0,170)
\put(0,0){\includegraphics[width=6cm,height=6cm]{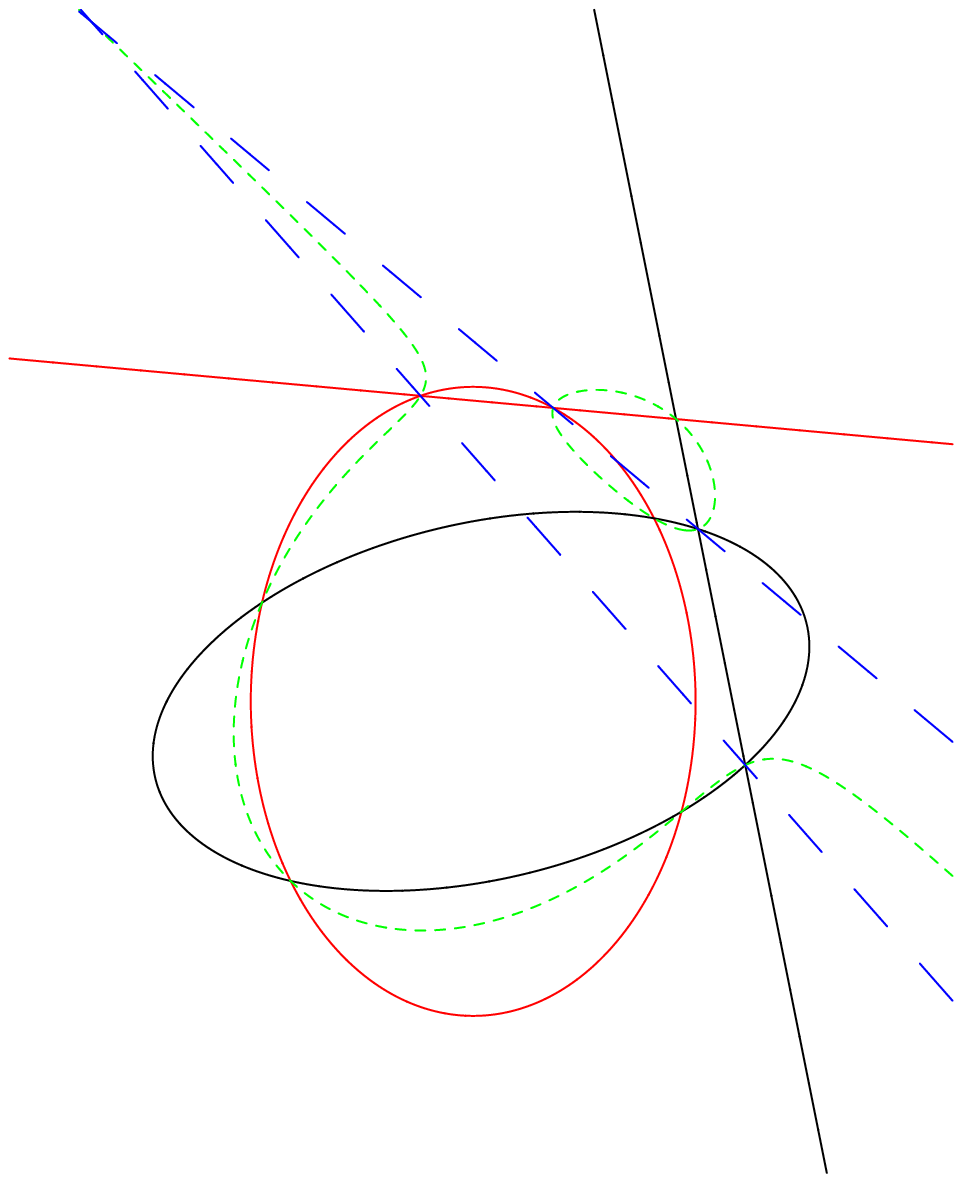}}
 \put(170,35){$E$}
 \put(125,115){{\tiny $p$}}
 \put(145,85){$\hat{Q}$}
 \put(80,15){$Q$}
 \put(125,47){{\tiny $\hat{q}_1$}} \put(130,95){{\tiny $\hat{q}_2$}}
 \put(95,117){{\tiny $q_2$}}
\put(10,160){{\tiny $o$}} \put(60,117){{\tiny $q_1$}}
\end{picture}

\begin{proof}

Recall that the intersection points $q_i$ of $Q \cap E$ are the
intersection points of the pairs of bitangents
$\theta_i,\theta_i+\ga$. With the notations of the previous
proof, if $p_1+p_2$ is the divisor associated to $\theta_i$, the
divisor associated to the product of the two lines is $2(p_1+p_2)+2
i(p_1+p_2)$. So $t_1=t_2$ and the points $q_i$ correspond to
singular points of singular conics of the form
$m((\lambda;\mu),(\lambda;\mu))$ i.e.
$$Q \cap E=
  \{(x;y;z)\; | \; \jac(A_1,A_2,A_3)^t \cdot(\lambda^2,\mu^2,2\lambda \mu)^t=0
\; \textrm{for one} \; (\lambda;\mu) \in \BP^1
  \}.
$$
Since $\frac{\partial A_2}{\partial x}=\frac{\partial A_1}{\partial y}=0$
the points of $Q\cap E$ are the $(x;y;z)$ coordinates for which  the
following system admits a solution
\[
  \begin{aligned}
  \frac{\partial A_2}{\partial x}\mu^2&+
  2\frac{\partial A_3}{\partial x}\lambda \mu=0,\quad
  &\frac{\partial A_1}{\partial y}\lambda^2+
  2\frac{\partial A_3}{\partial y}\lambda \mu=0\\
  \frac{\partial A_1}{\partial z}\lambda^2&+\frac{\partial A_2}{\partial z}\mu^2+
  2\frac{\partial A_3}{\partial z}\lambda \mu=0,\quad
  &(\lambda,\mu)\neq(0,0).
  \end{aligned}
\]
If $\lambda=0$ (resp. $\mu=0$) this system admits $q_1$ (resp.
$q_2$) as a solution. If $ \lambda \mu \ne 0$ then a solution
satisfies
\[
  (\frac{\partial A_2}{\partial x};-2\frac{\partial A_3}{\partial x})
  =(\lambda;\mu)
  =(-2\frac{\partial A_3}{\partial y};\frac{\partial A_1}{\partial y}),
\]
so it belongs to $\hat{Q}$. Hence the intersection points of $\hat{Q}$
and $Q$ are $E \cdot Q \setminus \{q_1,q_2\}=B'$ (see \ref{TagM}).\\
It follows from the definition of $\hat{q_i}$ that $\hat{q_i} \in \hat{Q}$
for $i=1,2$. We compute the intersections:
\[
  E \cdot (Q\ol{\hat{q_1}\hat{q_2}}) > (B' + q_1 +q_2) + \hat{q_1}+\hat{q_2}
\]
and
\[
  E \cdot (\hat{Q}\ol{q_1q_2}) > (B' + \hat{q_1} +\hat{q_2}) + (q_1+q_2+p).
\]
These $3$ cubics have thus $8$ points in common so
they lie in the same pencil. Moreover their last intersection point
is the same too, so $p \in \ol{\hat{q_1}\hat{q_2}}$.
\end{proof}

\begin{dsc}
{\em The curve $E'$.} For the purpose of describing the isotropic
subgroup $L_3$, as well as for technical reasons we set the
following notations: Denote by $Q_p$ the unique conic such that
$Q_p\cdot E=2p+B$, and by $Q'_p$ the unique conic such that
$Q'_p\cdot E'=2p+B'$ (the notations $B,B',p$ were defined in
\ref{Dlev} and \ref{Dcanon}). Recall that under the identification
of the linear system give by $k$, we have $E \cap Q \cap Q_p'=B'$ and
$E' \cap Q' \cap Q_p=B$.
\end{dsc}
\begin{proposition}\label{PcalctQ}
The plane cubics $E,Q T_p(E),\ol{q_1q_2} Q'_p$ lie on the same pencil.
\end{proposition}
 \hspace{3cm}
\begin{picture}(0,170)
\put(0,0){\includegraphics[width=6cm,clip]{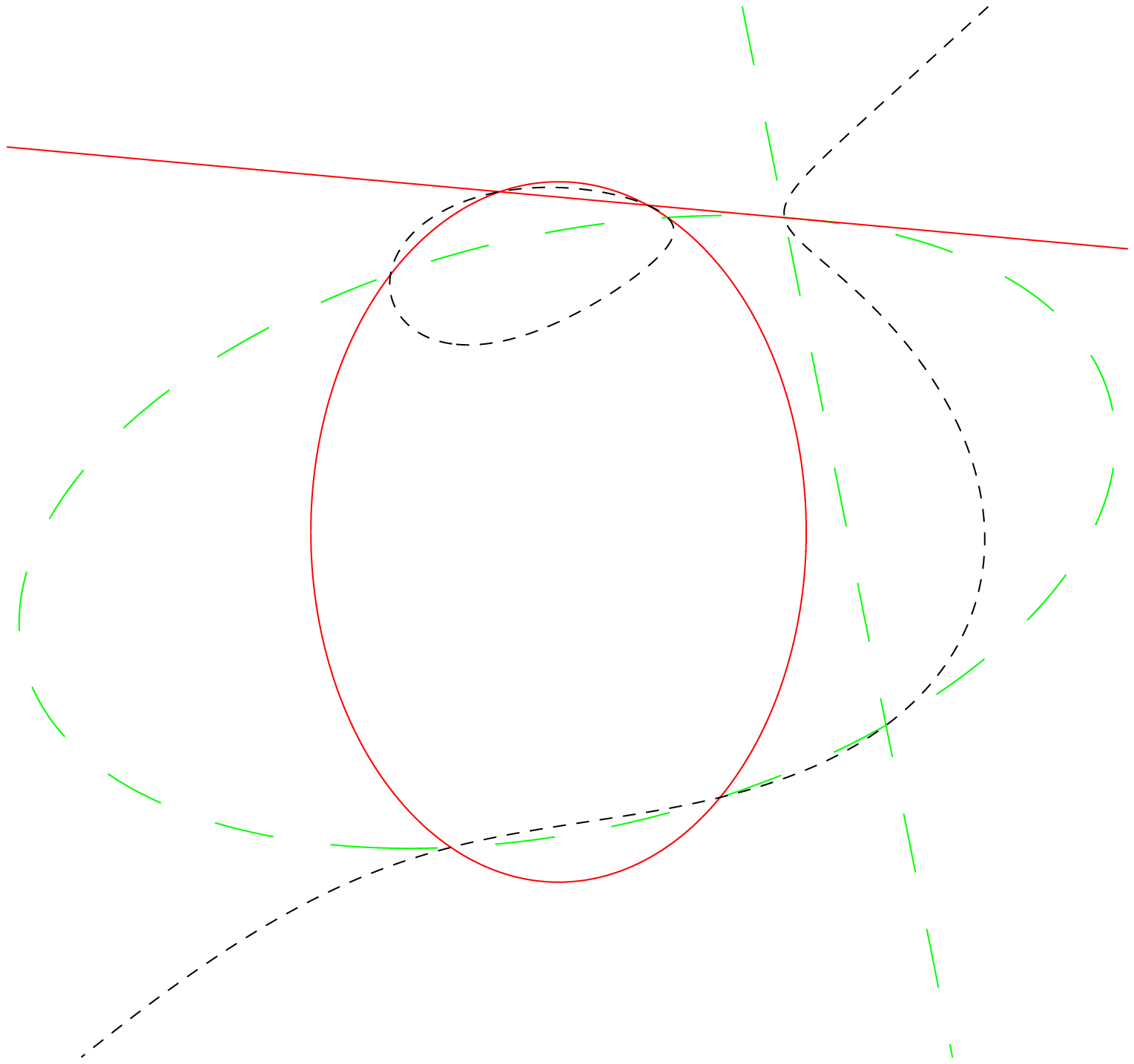}}
 \put(150,140){$E$}
 \put(125,127){{\tiny $p$}}
 \put(150,10){$T_p(E)$}
 \put(160,65){$Q_p'$}
 \put(80,15){$Q$}
\put(80,80){$B'$} \put(50,120){{\tiny $q_3$}} \put(60,25){{\tiny
$q_4$}} \put(105,30){{\tiny $q_5$}} \put(100,115){{\tiny $q_6$}}
\put(100,130){{\tiny $q_2$}} \put(66,135){{\tiny $q_1$}}
\end{picture}

\begin{proof}
As in the previous proof, this follows after calculating the intersections:
\[
  E\cdot(Q+T_p(E))> (B'+ q_1 + q_2) +2p,\quad
  E\cdot (\ol{q_1q_2}Q'_p) > (q_1+ q_2+p) + (B'+ p).
\]
\end{proof}
\begin{dsc} \label{PcalJ}
Let $J=\cup_{b\in B'}\ol{pb}$ ; we compute the quartic defining $J$ using the
following procedure.
\begin{enumerate}
\item Since $p=(-e;c;0)$, the lines passing through $p$ are given by
linear forms $c x+ e y - \ga z=0$ for some $\ga \in \BC$.
Thus $ J=\prod_{i=1}^4  (c x+ e y - \ga_i )$,
where the $\ga_i$'s are defined by the property $c x_i+ e y_i =\ga_i$
for each of the four points $(x_i;y_i;1) \in B'$.
\item let $Y=c x+e y$, then $Y_i:=c x_i+e y_i$ are the roots of
  the polynomial
  $R(Y)=(\textrm{Resultant}(Q(x,(Y-cx)/e,1),Q_p'(x,Y-cx)/e,1),x)$. So by
  definition, in {\em affine} coordinates, $J=R(c x+e y)$.
\end{enumerate}
\end{dsc}
\begin{proposition}\label{PcalcEQ}
The plane quartics $E'T_p(Q'_p),{Q'_p}^2,J$ lie on the same pencil.
\end{proposition}
 \hspace{3cm}
\begin{picture}(0,150)
\put(0,0){\includegraphics[width=6cm,clip]{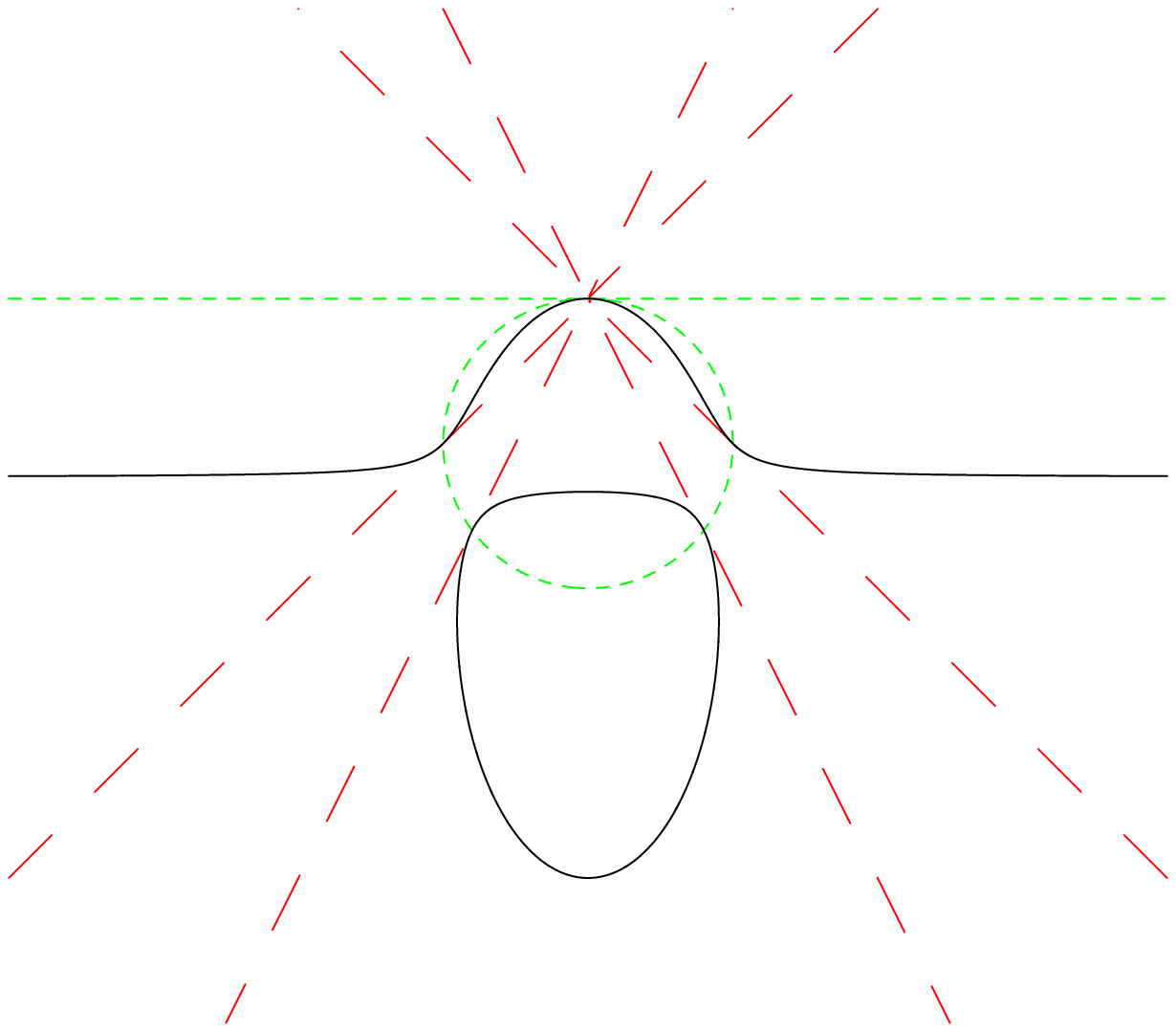}}
 \put(50,120){$J$}
 \put(97,107){{\tiny $p$}}
 \put(10,85){$E'$}
 \put(80,50){$Q_p'$}
\put(80,80){$B'$} \put(50,85){{\tiny $q_3$}} \put(60,70){{\tiny
$q_4$}} \put(105,70){{\tiny $q_5$}} \put(115,85){{\tiny $q_6$}}
\end{picture}

\begin{proof}\label{PcalcQQ}
By the definition of $J,Q'_p$ we have:
\[
 {Q'_p}^2\cdot J=2(B'+4p),\quad {Q'_p}^2\cdot (E'+T_p(Q'_p))=2((B'+2p)+2p).
\]
So the three quartics belongs to the same pencil.
\end{proof}

\begin{dsc} \label{Data14} {\em The curve $Q'$.} We compute $Q'$ using
Proposition \ref{PcalctQ} and symmetry. To compute the conic $Q_p$,
one notes that $Q_p \cdot E =2 p+B$. By definition of $B$ in
\ref{Dlev}, $b \in B$ if and only if $T_b(E) \cdot E=2b+p$. It is
then classical (see for instance Salmon \cite{salmon}) that $Q_p$ is
the polar conic of $p$ with respect to $E$ - recall that if
$p=(x_0;y_0;z_0) \in E$ then the polar conic of $p$ with respect to
$E$ is given by the equation $x_0 E_x + y_0 E_y + z_0 E_z=0.$
\end{dsc}

\begin{dsc}\label{DcalcAis}
To complete the calculation of $(C',\caL')$ from $(C,\caL)$ we still
have to make a final choice - the partition of the set $B'$ to the
two pairs $\{\{q_3,q_4\},\{q_5,q_6\}\}$. Geometrically, this is the
choice of the singular conic $\ol{q_3q_4}\cup\ol{q_5q_6}$ among the
three singular conics in the pencil spanned by the conics $Q,Q'_p$.
We start with a lemma on the symmetric situation:
\end{dsc}
\begin{lemma}\label{pointo}
The following equality holds:
\[
  \ol{q'_3 q'_4}\cap\ol{q'_5 q'_6}=p+\ga_E=o,
\]
where the addition is in $\Pic(E)$.
\end{lemma}
\begin{proof}
To see the first equality,  note that $2q'_i+p$ is a line section of $E\subset|K_C|^*$ for all
$i=3,\ldots,6$. By the proof of Theorem \ref{TAGM}, and the symmetry on the
construction, we also have $q'_3-q'_4=q'_5-q'_6=\ga_E$. Setting $\ti{o}:=p+\ga_E$ we see that
$q_3+q_4+\ti{o}=2q_3+p$ and $q_5+q_6+\ti{o}=2q_5+p$ are both line sections.
To see the second equality, note that $T_pE\cap T_o E$ lies on $E$, which means
that $\gc:=p-o$ is in $\Pic(E)[2]$. However, by a monodromy argument on maximal
isototropic flags on $\jac(C')[2]$ containing $L'_2$, we have $\gc=\ga_E$.
\end{proof}
\begin{dsc}
In order to find an equation for $(C',\ga')$ in the form $A_3'^2=A_1' A_2'$,
we apply projective transformation to mimic the form of the couple $(E,Q)$.
Denote by $T$
a projective transformation
of $|K_C|^*$ that sends $q'_1,q'_2$ and
$\ol{q_3q_4}\cap\ol{q_5q_6}$ to the points
$(1;0;0),(0;1;0)$ and $(0;0;1)$ respectively. Denote
$T(E')=(e_1,\ldots,e_{10})$ and $T(Q')=(d_1,\ldots,d_6)$.
Finally define $T_2$ to be the transformation which operates by a
multiplication of the
$x$-axis by $\sqrt{e_2/e_9}$, and a multiplication of the $y$-axis by
$\sqrt{e_4/e_6}$.
\end{dsc}
\begin{theorem}\label{TcalcAis}
The coordinates of the quadrics forms $T(A'_1),T(A'_2)$ are given by
\[
  (0,0,0,e_4/e_6,0,-1) \quad (e_2/e_9,0,0,0,0,-1),\\
\]
while the coordinates of the quadric form $T(A'_3)$ are given by
\[
   \frac{\begin{array}{l}
   (e_2 e_6(e_2 d_3-e_3 d_2),2e_2 e_3 e_4 d_6,2e_2^2 e_6 d_6,
    e_4(e_2 e_6 d_5-e_3 e_9 d_2)\\
  ,2 e_2 e_4 e_6 d_6,-e_6(- e_2 e_5 d_6+e_2 e_6 d_5+e_2 e_9 d_3-2 e_3 e_9 d_2)
   \end{array}}
{\sqrt{\begin{array}{l}
4 e_6 e_9 (e_6 e_2 e_5 d_6 e_3 d_2+e_6^2 e_2^2 d_3 d_5-
e_6^2 e_2 e_3 d_2 d_5 -e_9 e_6 e_2 d_3 e_3 d_2\\
+ e_9 e_6 e_3^2 d_2^2 -e_6 e_2^2 e_3 d_5 d_6 - e_6^2 e_2^2 d_2 d_6 +
e_2 e_3^2 e_4 d_6^2- e_2 e_3 e_4 e_6 d_3 d_6)
\end{array}}}
\]
\end{theorem}
\begin{proof}
Let us assume that we have taken $A_1'$ (resp. $A_2'$) such that
$T_2 \circ T(A_1')$ (resp. $T_2 \circ T(A_1')$) is the conic
$y^2-z^2$ (resp. $x^2-z^2$). Let $T_2 \circ
T(A'_3)=(a',b',c',d',e',f')$. By Theorem \ref{Tcoor}, if we call
$(E_2,Q_2)$ the data $(E,Q)$ associated to these transformations of
$A_1',A_2',A_3'$, we have
\[
  \begin{aligned}
    E_2&=(0;c';b';0;e';2(a'+d'+f');e';0;b';c';0),\\
    Q_2&=(0;c'e'({b'}^2-1+4a'd')-2b'({c'}^2d'+{e'}^2a');\\
     &b'(-2b'c'd'-e'+{b'}^2e')+2a'({b'}^2c'+c'{e'}^2-2b'e'f')-4{a'}^2b'e';0;\\
     &{b'}^3c'+2{c'}^2d'e'+2{b'}^2(d'-a')e'-b'c'(1+4{d'}^2+4d'f');\\
     &{c'}^2{e'}^2+{b'}^2({c'}^2+{e'}^2)-2b'c'e'(a'+d'+f')).
  \end{aligned}
\]
If we
let
 $T_2 \circ T(E')=(0,\hat{c},\hat{b},\hat{e},\hat{g},\hat{e},0,\hat{b},\hat{c},0)$
and $T_2 \circ T(Q')=(0,\delta_1,\delta_2,0,\delta_3,\delta_4)$,
 there exists a constant $\xi$ such that
\[
  b'=\hat{b}/\xi,\quad c'=\hat{c}/\xi,\quad e'=\hat{e}/\xi,
  \quad a'+d'+f'=\hat{g}/(2 \xi)
\]
and
\[
  \begin{aligned}
  \ &\left(\begin{array}{ccc}
   -2 \hat{b}{\hat{c}}^2 & -2 \hat{b}{\hat{e}}^2 & \hat{c}\hat{e}\\
   -2 {\hat{b}}^2\hat{c} & 2 {\hat{b}}^2\hat{c}+2 c{\hat{e}}^2-2\hat{b}\hat{e}\hat{g} & \hat{b}\hat{e} \\
   2 {\hat{b}}^2\hat{e}+2 {\hat{c}}^2\hat{e}-2 \hat{b}\hat{c}\hat{g} & -2 {\hat{b}}^2\hat{e} & \hat{b}\hat{c}
   \end{array}\right)
   \left(\begin{array}{c}d'\xi\\a'\xi\\4a'd'\xi^2-\xi^2\end{array}\right) +
   \\
     {\hat{b}}^2\left(\begin{array}{c}\hat{c}\hat{e}\\\hat{b}\hat{e}\\\hat{b}\hat{c}\end{array}\right)
    &=
    \frac{{\hat{c}}^2{\hat{e}}^2+{\hat{b}}^2({\hat{c}}^2+{\hat{e}}^2)-\hat{b}\hat{c}\hat{e}\hat{g}}{\gd_4}
    \left(\begin{array}{c}\gd_1\\ \gd_2\\ \gd_3\end{array}\right).
  \end{aligned}
\]
By the geometry of the configuration and the coordinates we chose, the
only solutions $a',d',\xi$ to the system above arise from solutions of the quadric
$T_2 \circ T(A'_3)$.
Since the quadric $A'_3$ is determined up to a sign, the matrix equation
above has only one solution, and this solution determines $a',d',\xi$ up to a
choice
of a sign. We apply then the transformation $T_2^{-1}$ to find the expression
of $T(A_3)$ in terms of $e_i,d_i$.
\end{proof}
\begin{remark}
Note that the transformation $T_2$, and the square roots
$\sqrt{e_2/e_9}$, $\sqrt{e_4/e_6}$ served merely as technical aids
in the proof above, and indeed vanished in the final result. The
situation is different with the root we take to distinguish between
$q'_1,q'_2$. Recall that when performing the trigonal construction,
one has to take a degree $2$ field extension in order to construct
$W'$ from the tower $Z'/X'/\BP^1$, and one has to construct $W'$ in
order to construct $C'$. Since after distinguishing between
$q'_1,q'_2$ we can construct $C'$, the root we take when we
distinguish between these points generate the field extension of the
function field of $W'$ over the function field of $Z'$. This carries
little significance when working over an algebraically closed field,
but when working over a non-algebraically closed field, it reflects
the fact that in order to find the isogeny $\prym(Z/X)\cong\jac(C)$
we may have to make a degree $2$ field extension of the base field.
\end{remark}
\begin{corollary}\label{Cdiff}
Let $\cM$ be the moduli of $a,b,c,d,e,f$ and a root of the cubic form (in $t$)
$\mathrm{Hessian}(t Q'_p+Q)$, then:
\begin{itemize}
\item The space $\cM$ is birational to a finite cover of the moduli of $C,\caL$
with monodromy group naturally isomorphic to $D_4$.
\item The map $T_2\circ T$ is defined globally over $\cM$. Moreover as a map
on quadrics in $x,y,z$ with parameters in $\cM$, the map
$T_2\circ T$ is an involution that lifts the involution $(C,\caL)\ra(C',\caL')$.
\item Using affine coordinates on $|K_{C}|^*$ (by setting $z=1$) the map $\ol{k}$ is given by the formula
\[
  T_2\circ T\left(\frac{l dx}{\partial(A'_1A'_2-{A'_3}^2)/\partial y}\right)
  =\pm{\ol{k}}^{-1}\left(
  \frac{l dx}{\partial(A_1A_2-A_3^2)/\partial y}\right),
\]
where $l$ is any linear form.
\end{itemize}
\end{corollary}
\begin{proof}
The first assertion follows from the choice of coordinates we use (see
Theorem \ref{Tcoor}), and the fact that the singular conics in
the pencil of conics spanned by $Q'_p,Q$ are in 1-1 correspondence with
the roots of the cubic $\mathrm{Hessian}(t Q'_p+Q)$. The dihedral group is
the symmetry group acting on the nested partition of linear forms
$\{\{x-z,x+z\},\{y-z,y+z\}\}$.\\
The second assertion follows from
the definition of $T$ and $T_2$, and from Theorem \ref{TcalcAis}.\\
It is well known that a basis of regular differentials on a genus
$3$ non hyperelliptic curve $C$ can be given by $\left(\frac{l
dx}{\partial(C)/\partial y}\right)$. With the identifications we
have made during the construction on the coordinates
$(x;y;z),(x';y';z')$ (see \ref{Dprog}), the map $\ol{k}$ with this
choice of bases is given by the transformation $(T_2 \circ T)^{-1}$
up to a constant. However, since $T_2\circ T$ is an involution on
$\cM$, the square of this constant is $1$.
\end{proof}
\begin{dsc}\label{Dcalcbasis}
Our final objective in this section is to show how one iterates the
construction. Following \ref{Dbasis}, we denote with
$\tilde{\phantom{1}}$ the objects related to $(C',\ti{\caL})$.\\
Our first task is to find an $\ti{\ga}$. By our analysis of the symplectic
pairings in Proposition \ref{Plev} and in \ref{Dlev} we have
\[
  \begin{aligned}
  \# D^{-1}_{\ga'}(\ti{\ga}) \cap \{q_1',q_2'\} = 0 & \iff&
         \ti{\ga} \in L_2'^\perp \sm L_2', \\
  \# D^{-1}_{\ga'}(\ti{\ga})
     \cap \{q_5',q_6'\} =1 & \iff & \ti{\ga} \notin L_3'.
  \end{aligned}
\]
So we can assume that
$D^{-1}_{\ga'}(\ti{\ga})=\{q_3',q_5'\}$. The situation can be
represented as follows :
\begin{figure}[H]
\begin{center}
\label{fig}
\input{figtilde.pstex_t}
\end{center}
\end{figure}
Thus, the four bitangents $\gb'_{31},\gb'_{32}$ (lying over $q'_3$)
$\gb'_{51},\gb'_{52}$ (lying over $q'_5$) in
$\gC_{\ga'\oplus\ga}=\gS_{\ga'} \cap \gS_{\ti{\gs}}$ (Proposition \ref{P46})
can be grouped as
$(\gb'_{31},\gb'_{52})$ and $(\gb'_{51},\gb'_{32})$ to give
a possible $\ti{\ga}$ (the other grouping correspond to an
$\ti{\ga}+\ga'$).
Let us denote $\ti{q}_3=\gb'_{31} \cap
\gb'_{52}$ and $\ti{q}_5=\gb'_{51} \cap \gb'_{32}$. Let
$\ti{E},\ti{Q}$ be the cubic and conic associated to
$(C',\ti{\ga})$. We denote $\{\ti{q}_i\}_{i=1}^6=\ti{Q} \cap \ti{E}$.

The second step is to find $\ti{L}_2$.
By Proposition \ref{Plev}, a maximal isotropic group which contains
$\ga'$ (resp. $\ti{\ga}$) is equivalent to the partition of
$\gS_{\ga'}/\ga'$ (resp.
$\gS_{\ti{\ga}}/\ti{\ga}$) in three pairs. There are three
different maximal isotropic spaces containing $\ga'\oplus\ti{\ga}$, given by
the non-zero points in $(\ga' \oplus \ti{\ga})^{\perp}/(\ga' \oplus \ti{\ga})$.
These maximal flags are
in bijection with partitions two pairs of the points $\{q'_1,q'_2,q'_4,q'_6\}$,
and also with partition to pairs of the points
$\{\ti{q}_1,\ti{q}_2,\ti{q}_4,\ti{q}_6\}$.
One of the three maximal isotropic groups containing $\ti{\ga}\oplus\ga'$ is
$\ti{\ga} \oplus L_2'$, the two others corresponds to
$\ga'\oplus\ti{L}_2$ for the two different choices of $\ti{L}_2$.
Thus, in order to choose $\ti{L}_2$ we first choose the maximal isotropic group
$\ga'\oplus\ti{L}_2$, and in making this choice, we exclude the partition
corresponding to the group $\ti{\ga} \oplus L_2'$.\\
The partition of $\{q'_1,q'_2,q'_4,q'_6\}$ which correspond
to $\ga'\oplus\ti{L}_2$ is simply the partition $\{q'_1,q'_2\},\{q'_4,q'_6\}$,
but in order to proceed we will have to find the corresponding partition
of $\{\ti{q}_1,\ti{q}_2,\ti{q}_4,\ti{q}_6\}$. To do this we will describe
explicitly the natural isomorphism between the three partitions to two pairs of
these 4-tuples. This isomorphism is geometric in nature, and to describe
it we will interpret these partitions as singular conics defined by the
partitions\\
To describe the isomorphism we need some more notations : if
$\ga_1,\ga_2 \in \jac(C')[2]$ such that $\langle
\ga_1,\ga_2 \rangle =0$, one denote by $A_{3,\ga_i},
E_{\ga_i},Q_{\ga_i},p(\ga_i)$ the elements $A_3,E,Q,p$
relative to the construction starting from $(C',\ga_i)$. Also we
denote by $E_{\ga_i \oplus \ga_j}$ and
$Q_{p(\ga_i),\ga_j}$ the curve $E'$ and the conic $Q_p'$
constructed from the data $(C',\ga_i \subset \ga_i \oplus
\ga_j)$ and by $'$ the symmetric constructions. Note that
$$E_{\ga_i \oplus \ga_j} = E'_{(\ga_i \oplus
\ga_j)^\perp/L} = E_{\ga_j \oplus \ga_i}$$ for any maximal
isotropic group $L$ containing $\ga_i \oplus \ga_j$. However,
$Q_{p(\ga_i),\ga_j} \ne Q_{p(\ga_j),\ga_i}$.
\end{dsc}
\begin{lemma}
The identification  of a singular conic through
$\{q'_1,q'_2,q'_4,q'_6\}$
 (resp. $\{\ti{q}_1,\ti{q}_2,$ $\ti{q}_4,\ti{q}_6\}$) with a root of
$$\textrm{Hessian}(Q_{\ga'} + u Q_{p(\ga'),\ti{\ga}})
\; \textrm{(resp.} \; \textrm{Hessian}(Q_{\ti{\ga}} +u
Q_{p(\ti{\ga}),\ga'}) \textrm{)}$$ defines a natural
transformation $\mu : \BP^1 \to \BP^1$ which fixes $\infty$ and maps
the two triple of roots bijectively.
\end{lemma}
\begin{proof}
We have $E_{\ga' \oplus \ti{\ga}}=E_{\ti{\ga} \oplus\ga'}$. Let $t$ be the
translation on $E_{\ti{\ga} \oplus\ga'}$ by $p(\ti{\ga})-p(\ga')$. Let
$s_{Q_{p(\ga'),\ti{\ga}}},s_{Q_{\ga'}}$ (resp.
$s_{Q_{p(\ti{\ga}),\ga'}},s_{Q_{\ti{\ga}}}$) denote the
sections of the bundle of $E_{\ga' \oplus \ti{\ga}}$ defined
by the divisor $q_1'+q_2'+q_4'+q_6'$ (resp.
$\ti{q}_1,\ti{q}_2,\ti{q}_4,\ti{q}_6$) and corresponding to the
subscript objects. Since  $Q_{p(\ga'),\ti{\ga}}$ (resp.
$Q_{p(\ti{\ga}),\ga'}$) is \emph{the} polar conic of
$p(\ga')$ (resp. $p(\ti{\ga})$), there exists $a \in K$ such
that
\[
 t^*(s_{Q_{p(\ga'),\ti{\ga}}})=a s_{Q_{p(\ti{\ga}),\ga'}}.
\]
Thus $t$ maps the points
$\{q'_1,q'_2,q'_4,q'_6\}$ onto  $\{\ti{q}_1,\ti{q}_2,$
$\ti{q}_4,\ti{q}_6\}$. Since $t^*(s_{Q_{\ga'}})$ contains the
points $\{\ti{q}_1,\ti{q}_2,\ti{q}_4,\ti{q}_6\}$, there are also two
constants $b,c$ such that $t^*(s_{Q_{\ga'}})=b
s_{Q_{\ti{\ga}}}+ c s_{Q_{p(\ti{\ga}),\ga'}}.$ Thus any
section in the pencil $s_{Q_{\ga'}}+u
s_{Q_{p(\ga'),\ti{\ga}}}$ is map through $t$ on $b
s_{Q_{\ti{\ga}}}+ (au+c) s_{Q_{p(\ti{\ga}),\ga'}}$. Hence
there is an affine transformation $\mu$ which maps a conic
$Q_{\ga'} + u Q_{p(\ga'),\ti{\ga}}$ to $Q_{\ti{\ga}'} +
\mu(u) Q_{p(\ti{\ga}),\ga'}$. In particular, a singular conic
is mapped to a singular conic which means that $\mu$ maps the three
roots of $\textrm{Hessian}(Q_{\ga'} + u
Q_{p(\ga'),\ti{\ga}})$ to the three roots of
$\textrm{Hessian}(Q_{\ti{\ga}} + u
Q_{p(\ti{\ga}),\ga'})$.\\
To identify the transformation $\mu$, we work on the generic case
$C'$ given by $(a',b',c',d',e',f')$ and we are looking for a
continuous affine transformation. We assume that the bitangents
$\gb'_{31},\gb'_{32},\gb'_{51},\gb'_{52}$ are
$y-z,y+z,x-z,x+z$ respectively.
 Denote
by
\[
  T:=\frac{1}{2}
 \left(\begin{array}{rrr}-1&1&2\\1&-1&2\\1&1&0\end{array}\right)
\]
the transformation which sends $y-z,y+z,x-z,x+z$ to
$y-z,x+z,x-z,y+z$ respectively. This projective transformation
defines a linear transformation $T_{\ti{\ga}}$ on the coefficients
of $A_{3,\ga'}$ given by $(a',b',c',d',e',f')$ maps to the
coefficients of $A_{3,\ti{\ga}}$
\[
  \begin{aligned}
  (\frac{1}{4}&(a'-b'-c'+d'+e'+f'),\frac{1}{2}(-a'+b'-d'+f'),
  \frac{1}{2}(-2a'+c'+2d'+e'),\\
   \frac{1}{4}&(a'-b'+c'+d'-e'-f'),\frac{1}{2}(2a'+c'-2d'+e'),a'+b'+d')
  \end{aligned}.
\]
Using Theorem \ref{Tcoor}, we can compute the different objects
involved and we find \[
  \begin{aligned}
  \mu(u)=&(2a'-c'-2d'-e')\cdot (2a'+c'-2d'+e')  \\
    &\cdot (2a'b'{e'}^2-4a'c'd'e'-{b'}^2c'e'+2b'{c'}^2d'+c'e') u\\
     &+2 e' \cdot c' \cdot a' \cdot \\
  &\cdot   (4{a'}^2b'-8a'b'd'-2a'c'e'-b'{c'}^2+4b'{d'}^2-b{e'}^2-2c'd'e'+2c'e'f').
  \end{aligned}
\]
Let us say a few words about the computation: we are comparing the
coefficients of two monic cubic forms in $u$, under the
transformation $u \mapsto (\mu_0 u+\mu_1)$. Thus we get the
equation:
\[
  \epsilon (u^3+a_2 u^2+a_1 u+a_0)=(\mu_0 u+\mu_1)^3+
  b_2 (\mu_0 u+\mu_1)^2+b_1 (\mu_0 u+\mu_1)+ b_0.
\]
Comparing the $u$ coefficients we get a system of equations in $\mu_0,\mu_1$:
$$\begin{cases}
3 \mu_1 + b_2- \mu_0 a_2 &=0, \\
3 \mu_1^2 + 2 b_2 \mu_1+ b_1- \mu_0^2 a_1 &=0,\\
\mu_1^3 + b_2 \mu_1^2 + b_1 \mu + b_0- \mu_0^3 a_0 &=0.
\end{cases}$$
We solve the system by finding the two solutions of the first two equations,
and checking which of the two solutions solves the third equation.
\end{proof}
\begin{dsc} \label{Dexpbasis}
By a projective transformation one can send the bitangents
$\gb'_{31},\gb'_{32}$, $\gb'_{51},\gb'_{52}$ to
$y-z,y+z,x-z,x+z$ respectively. Using the previous lemma, one can
then identify the value $u_0$ of $u$ corresponding to
$\ti{\ga}+L_2'$ (i.e to the singular conic whose one component is
$\ol{q_1'q_2'}$) and then exclude the singular conic
$Q_{\ti{\ga}} + \mu(u_0) Q_{p(\ti{\ga}),\ga'}$. Let us
denote this one $\ol{\ti{q}_1\ti{q}_2} \cup
\ol{\ti{q}_4\ti{q}_6}$ and then $\ti{L}_2$ is represented for instance
by $\ol{\ti{q}_4\ti{q}_6}$.\\
 The last task is to identify
 $\ti{L}_3$ : by Lemma \ref{L2t}, it has
 to contain one of the points $\ti{q}_3,\ti{q}_5$ so it is given by any
 of the choice of a
 pair $\{\ti{q}_2,\ti{q}_5\}$ or $\{\ti{q}_2,\ti{q}_3\}$.
\end{dsc}
%
\section{Real curves}\label{Sreal}
%
\begin{dsc}
In this section we show that if $C$ is a real $M$ curve of genus $3$
(i.e. a curve with $4$ components), then the {\em topology} of the
real structure induces
a distinguished isotropic flag $L_1\subset L_3$ in $\jac(C)[2]$ such
that the curve $C'$ is an $M$ curve. In Theorem \ref{TRagM} we
establish a bijection between partitions of the four components of
the curve $C$ to two pairs, and the set of full flags $\caL$
containing the flag $L_1\subset L_3$. In \ref{Diter1} we show how to
find the {\em topologically} distinguished flag $\ti{\caL}'$ on the curve $C'$
using the data $C',\caL'$ (one calculates $C',\caL'$ from the pair $C,\caL$)
by taking square roots as described in \ref{Dcalcbasis}
 - thus getting an iterative process - we will show that the choice of these
square roots is uniquely determined by the topology.
Finally, in \ref{Drint} we describe the iterative integration algorithm.
\end{dsc}
\begin{dsc}\label{DMcurve}
Let $C$ be a real plane quartic with $4$ components $C_1,C_2,C_3,C_4$.
We denote by $\jac_\BR(C)$ the real part of the Jacobian of
the curve $C$, and by $\jac_\BR(C)_0$ the $0$-component of  $\jac_\BR(C)$ -
see \cite{GH} p. 159.
Recall that since the degree of the curve $C$ is even, each of the $C_i$s
is null homotopic in $\BP\BR^2$ (see e.g. \cite{GH}).
Whence, the set $\BR\BP^2\sm C_i$ is a union
of a disk and a M\"obius strip. Recall also that
the quotient $\jac_\BR(C)[2]/\jac_\BR(C)_0[2]$ is naturally isomorphic to the vector
space $\BF_2[C_1,C_2,C_3,C_4]/\BF_2$,
where $\BF_2$ acts by adding $1_{\BF_2}$ to all the coordinates.

Let $\{\{C_1,C_3\},\{C_2,C_4\}\}$ be a partition of the four
components to two pairs. Denote by $c_i$ a point in the trivial
component of $\BP\BR^2\sm C_i$, and choose the infinity line
$l_\infty$ in $\BP\BR^2$ such that the $4$ points $c_1,c_2,c_3,c_4$
admit a cyclic order in $\BP\BR^2\sm l_\infty$ (formally, this means that $c_1,c_2,c_3,c_4$ sit on an ellipse in $\BP\BR^2\sm l_\infty$ in the order $1,2,3,4$).
We assume the order induced on $c_1,c_2,c_3,c_4$ from the choice of the
line $l_\infty$ is counter clock wise. Note that there is a natural
isomorphism
\[
  \mathcal{H}:=H_1(\BP\BR^2\sm\{c_1,c_2,c_3,c_4\},\BF_2)\cong
  \BF_2[l_\infty]\oplus\BF_2[C_1,C_2,C_3,C_4]/\BF_2.
\]
Since the bitangents are lines, the classes of the bitangents in $\mathcal{H}$
have non trivial $l_\infty$ coordinate; in fact a much stronger result holds:
\end{dsc}
\begin{theorem}\label{TRagM}
For any $i\in\{1,\ldots,4\}$ there is a Steiner system $\gS_i$ such that
the bitangents in $\gS_i$ have exactly four representatives in each of the
following homology classes in $\mathcal{H}$:
\[
  l_\infty+C_i,\quad l_\infty+C_i+C_{i+1 \pmod{ 4}}, l_\infty+C_{i-1 \pmod{4}}.
\]
\end{theorem}
\begin{proof}
All additions of indices in the proof are modulo $4$. For any $i<j$
consider a one parameter degeneration of the curve $C$ to a curve
$\ol{C}_{ij}$ such that on $\ol{C}_{ij}$, the ovals $C_i,C_{i-1}$
are connected with a node, and the ovals $C_j,C_{j-1}$ are connected
with a node - see the figure below.
\begin{center}
\includegraphics[angle=0, width=0.5\textwidth]{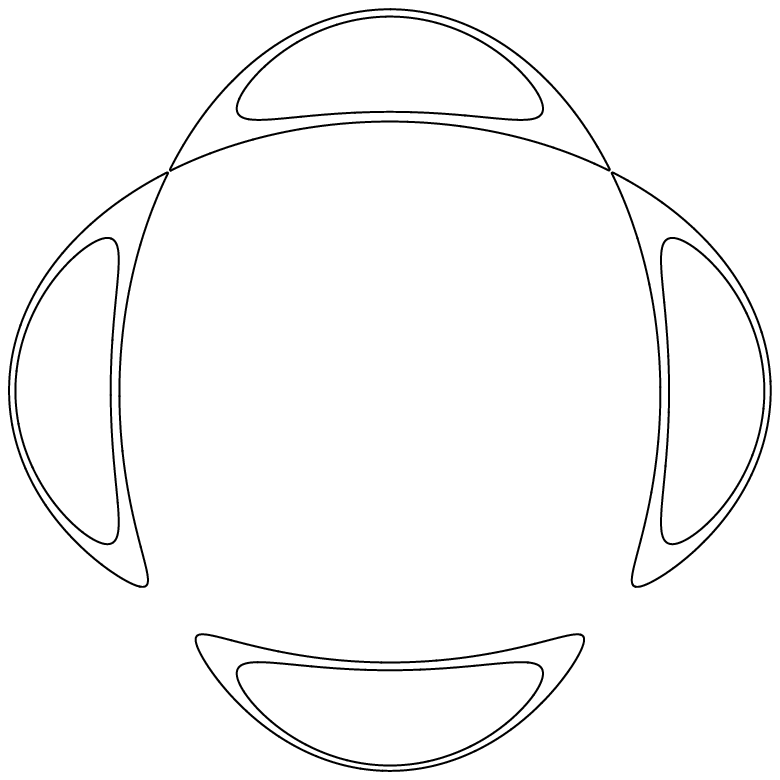}
\end{center}
By \cite{Be1} the degeneration of  each of the two nodes degenerates
a $2$-torsion point in the Jacobian variety $\jac(C)$; we mark these
points by $\gc_i,\gc_j$. Under the degeneration from $C$ to
$\ol{C}_{ij}$, any of the $6$ pairs of bitangents in the Steiner
system $\gS_{\gc_i}$ degenerates to one double line through the node
of corresponding to $\gc_i$; the same property holds also for $j$.
Moreover the intersection of the Steiner systems
$\gS_{\gc_i}\cap\gS_{\gc_j}$ is then the quadrapole line through the
two nodes of $\ol{C}_{ij}$. Thus by Proposition \ref{Plev} the Weil
pairing $\langle\gc_i,\gc_j\rangle$ is $0$, and all the bitangents
in $\gC_{\gc_i\oplus\gc_j}$ have the same homology class in
$\mathcal{H}$: If $j=i+1$ this class is $l_\infty+C_i$, and if
$j=i+2$ this class is $l_\infty+C_i+C_{i+1}$.
\end{proof}
\begin{dsc}
In Theorem \ref{TRagM} we identified the homology classes of $24$ bitangents.
Since for any $i$
there is a bitangent to $C_i,C_{i+1}$ in the class $l_\infty$, this completes
the identification of the bitangents. By Proposition \ref{Plev} we have also
identified a distinguished flag given by the following partition data:
\begin{itemize}
\item Four bitangents in the class of $l_\infty$.
\item Four bitangents in each of the classes $l_\infty,l_\infty+C_1+C_3,l_\infty+C_2+C_4$
- by the combinatorial structure we described and by Proposition
\ref{Plev} this is a Steiner system.
\item Four bitangents in each of the classes in $\mathcal{H}$ with non trivial line
coordinate except the class $l_\infty+C_1+C_3$.
\end{itemize}
Finally recall (see \cite{Hui}) that since $C$ is an $M$ curve, the
variety $\jac_\BR(C)_0$ is naturally isomorphic to a product of any
$3$($=\text{genus}(C)$) of the components. Since the group
$\jac_\BR(C)_0[2]$ is spanned by any two subgroups of order $4$ in
it, and since four of the order $4$ subgroups we built above sit in
the product of $3$ components of the curve, the maximal group in the
flag we built is $\jac_\BR(C)_0[2]$. Note that the quotient
$\jac_\BR(C)/L_3$ is the Jacobian of an $M$ curve if and only if
this quotient has $2^3$ components, which means that
$L_3=\jac_\BR(C)_0[2]$ - thus the choice of a distinguished $L_3$ we
made above is indeed the unique choice which will enable iteration.
\end{dsc}
\begin{dsc}\label{Diter1}
To complete the description of an iterative algorithm we have to solve
two problems: initiating the algorithm, and performing an iterative
step. In the discussion of these problems we will apply several times the
proposition below.
\end{dsc}
\begin{proposition}\label{Prealdef}
The symplecto-algebraic properties of a configuration of bitangents to a real
M-curve of genus $3$, as points on the odd part of an affine symplectic
space, are determined by
the homotopy classes of the bitangents in $\BP\BR^2\sm\{c_i\}_{i=2}^4$
and the intersection pattern of the bitangents with the components of $C$.
\end{proposition}
\begin{proof}
This follows from the following facts:
\begin{itemize}
\item Bitangents are continuous on families. Thus, the homology classes of
bitangents in $\BR\BP^2\sm\{c_i\}_1^4$ is constant on families.
\item Level structure is continuous on families.
\item The moduli space of $M$-curves is irreducible.
\end{itemize}
\end{proof}
\begin{dsc}
We will apply Proposition \ref{Prealdef} several times to
study the configuration of bitangents arising from
\ref{Dcalcbasis}. We will describe real algebro-geometric data on the moduli
of configurations of bitangents that defines
several non-zero real algebraic functions.
To show that some real configuration is associated with
a distinguished flag (in the sense of Theorem \ref{TRagM}) we will present one curve
for which our function is positive on the distinguished configuration and
negative on the others.
The conceptual calculations appear below. The related numeric calculations
are in \cite{LR}.

\item{\em Initiating the algorithm.}
To initiate the algorithm one essentially has to solve, in a Galois theory
sense, the bitangents of the curve $C$.
As the Galois group acting on the bitangents is generically unsolvable,
this problem is generically unsolvable in radicals. However, there are
still other computationally useful problems we will answer:
\item $\bullet$ Determine if on a given M-curve of genus $3$ and quadrics
$A_1,A_2,A_3$, the flag $L_1\subset L_2$ induced from $A_1,A_2,A_3$ is a
subflag of a distinguished flag (in the sense of Theorem \ref{TRagM}).
To check this it suffices
to verify that the curve $C$ lies in one component of
$\BP\BR^2\sm\mathrm{Nulls}(A_1A_2)$. It suffices to check this on an
infinitesimal neighborhood of the $4$ bitangents determined by the
conics $A_1,A_2$.
\item $\bullet$ Given an M-curve $C$ in the form $A_1A_2-A_3^2$ such that the
induced flag is a subflag of a distinguished flag, mark the choice of
the distinguished group $L_3$.
Note that the quartic form
$(x^2-z^2)(y^2-z^2)$ separate the real projective plane to $3$ positive
components and $4$ negative components. Thus, if the form $A'_3$,
calculated as in \ref{TcalcAis} is purely imaginary, then the curve $C'$
is an M curve, and thus the choice of $L_3$
is the distinguished choice.
Note also that $A_3'$ is
purely imaginary if and only if the expression  under the square root in
Theorem \ref{TcalcAis} is negative. By Proposition
\ref{Prealdef}, it suffice to show one example of a curve $C$ with purely
imaginary $A'_3$. We do this in \cite{LR}.
\item{\em Describing the iterative step:}
Recall (see \ref{Dcalcbasis}) that during the calculation of the
iterative step in \ref{Dcalcbasis}, the field extensions where
geometrically described by making several times choices of the
following type: given three lines $l_1,l_2,l_3$ and two pairs of
points $p_{i1},p_{i2}\in l_i\sm(l_i\cap(l_{3-i}\cup l_3))$ find a
partition to two pairs of the four points
$p_{11},p_{12},p_{21},p_{22}$ that is not the one arising from the
lines $l_1,l_2$. This choice boils down to a positivity question: we
consider the pencil of conics throughput the four points
$p_{11},p_{12},p_{21},p_{22}$. There are three singular conics in
this pencil, one of which is given by $l_1\cup l_2$. Moreover,
$l_1\cup l_2$ cuts $\BP\BR^2$ to two components, and the two nodes
of the two other singular conics in the pencil, which are
$\ol{p_{11}p_{21}}\cup\ol{p_{22}p_{12}},\ol{p_{11}p_{22}}\cup\ol{p_{12}p_{21}}$
appear one in each of the two components of $\BP\BR^2\sm(l_1\cup
l_2)$. We calculate the choices which bring us to a distinguished on
$C'$ (in the sense of Theorem \ref{TRagM} in \cite{LR}). We plot in
the figure below one step of the computation - finding the
distinguished topological configuration of
$\gC_{L'_1\oplus\ti{L}_1}$ (See \ref{Dcalcbasis}). The set
$\gC_{L'_2}$ is plotted in dashed red lines and the set
$\gC_{\ti{L}_1\oplus L'_1}$ is plotted in blue dotted lines.
\begin{center}
\begin{picture}(130,170)
\put(0,0){\includegraphics[width=6cm,clip]{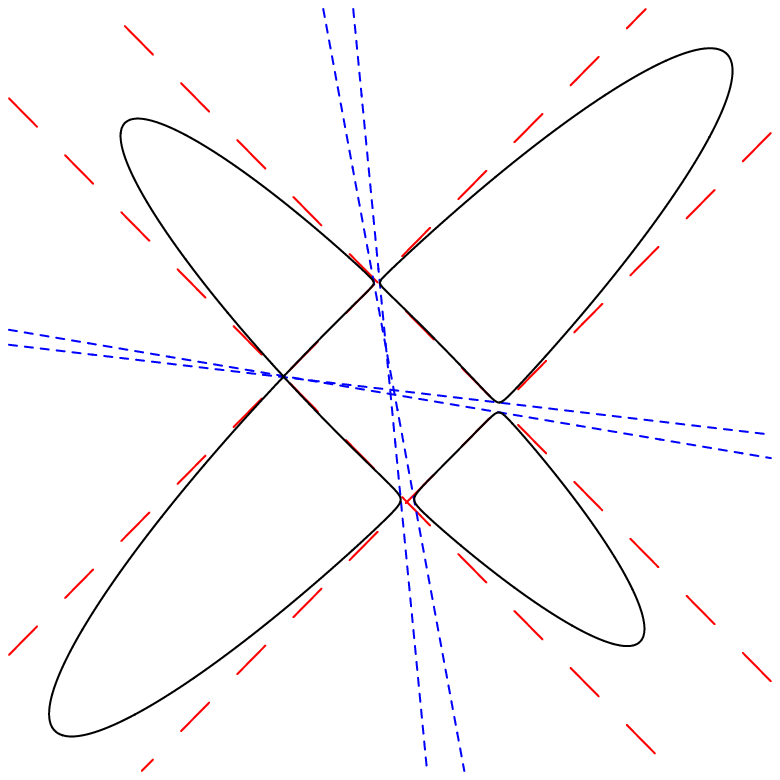}}
 \put(130,20){$\gC_{L_2}$}
 \put(130,85){$\gC_{\ga'\oplus\ti{\ga}}$}
\end{picture}
\end{center}
\end{dsc}
\begin{dsc}\label{Drint}
We conclude this section with a description of the integration algorithm.
By \cite{Hui} the $0$ component of the real Jacobian of the curve $C$ is
a product of any $3$ of the components of $C$. Thus, real projective
cycles on the Jacobian of $C$ are homologous to sums of $C_i$s.
By the same argument, these sums map under our construction to
sums of the components of $C'$.
Considering the theta function of $C$ we see that when the AGM algorithm
is applied iteratively, the distances between the
bitangents in each of the distinguished $4$-tuples are decreasing
exponentially.
i.e., the limit curve of this process is a union of $4$ lines in $\BP\BR^2$.
Using this method and Corollary \ref{Cdiff} we reduced the calculation of
integrals of cycles on
$\jac_{\BR}(C)$ to calculation of integrals of rational functions on line
segments.

In the figure below, which is computed in the final step of
\cite{LR}, we plot two
iterations of our algorithm, where the canonical classes of the
curves $C,C'$, and the next curve in the iterative process are
identified.
\begin{center}
\begin{picture}(130,170)
\put(0,0){\includegraphics[width=6cm,clip]{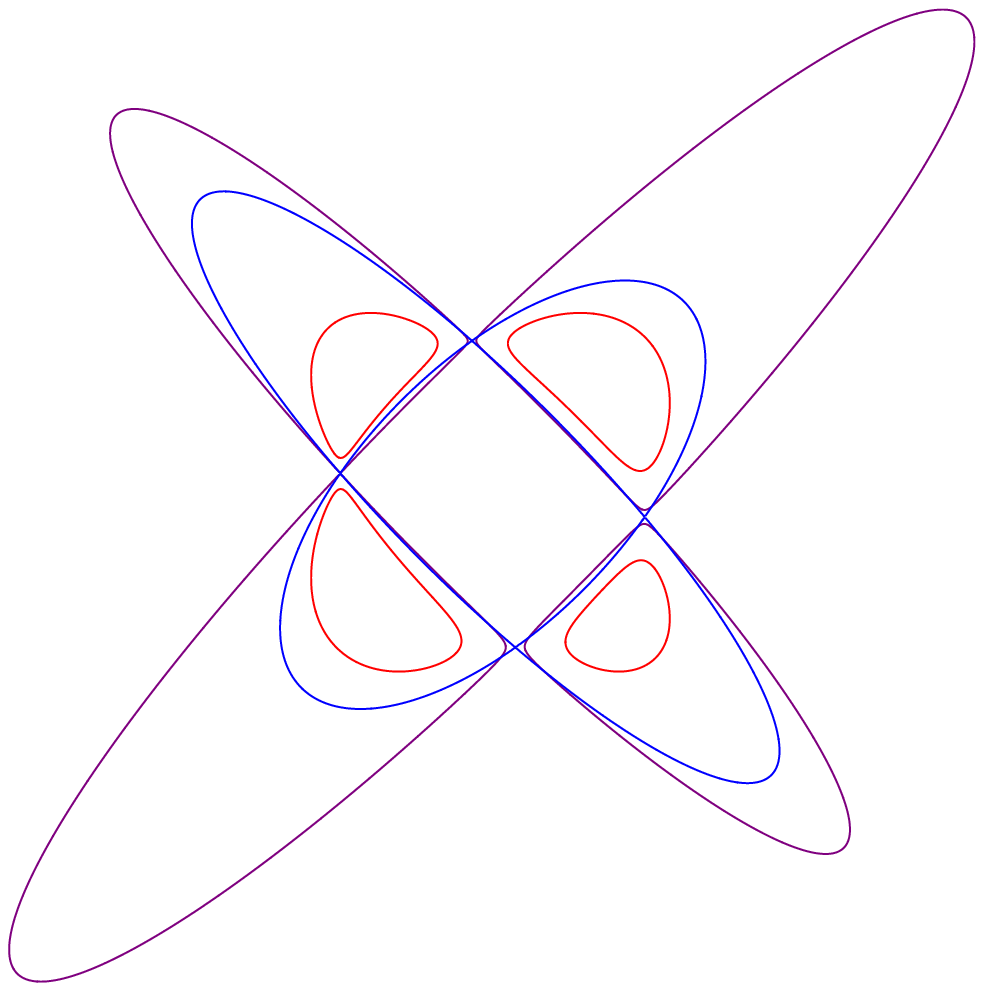}}
 \put(100,100){$C$}
 \put(145,145){$C'$}
 \put(120,120){Second iteration}
\end{picture}
\end{center}
\end{dsc}



\begin{thebibliography}{99}

\bibitem[ACGH]{ACGH}  Arbarello, E.; Cornalba, M.; Griffiths, P. A.; Harris, J. {\em Geometry of algebraic curves. Vol. I}. Grundlehren der Mathematischen Wissenschaften, 267. Springer-Verlag, 1985.
\bibitem[Bea1]{Be1} Beauville, A. {\em Prym varieties and the Schottky problem}. Invent. Math. 41 (1977), no. 2, 149-196.
\bibitem[Bea2]{Be2} Beauville, A. {\em Sous-vari\'et\'es sp\'eciales des vari\'et\'es de Prym}. Compositio Math. 45  (1982), no. 3, 357-383.
\bibitem[BB]{BB} Borwein, J. M.; P. B. Borwein, P. B. {\em Pi and the AGM: A Study in Analytic Number Theory and Computational Complexity}. A Wiley-Interscience Publication. John Wiley \& Sons, Inc., 1988.
\bibitem[BM]{BM}  Bost, J.-B.; Mestre, J.-F. {\em Moyenne Arithmetico-geometrique et P\'eriodes des Courbes de genere 1 et 2}.
Gaz. Math. No. 38 (1988), 36-64.
\bibitem[Cob]{Co} Coble, A. {\em Algebraic geometry and Theta functions}. Revised printing. American Mathematical Society Colloquium Publication, vol. X American Mathematical Society, 1961.
\bibitem[Cox]{Cox} Cox, D. A. {\em The arithmetic-geometric mean of Gauss}.  Enseign. Math. (2)  30  (1984),  no. 3-4, 275-330.
\bibitem[Dol1]{Dol1} Dolgachev, I. {\em Rationality of $\caR_3$}.
              Unpublished notes, available online at
              \href{http://www.math.lsa.umich.edu/~idolga/preprints.html}{http://www.math.lsa.umich.edu/\~{}idolga/preprints.html}.
\bibitem[Dol2]{Dol2} Dolgachev, I. {\em Classical Algebraic geometry}.
              A manuscript in progress. Available online at
              \href{http://www.math.lsa.umich.edu/~idolga}{http://www.math.lsa.umich.edu/\~{}idolga}.
\bibitem[Don]{Don}  Donagi, R. {\em The fibers of the prym map}.
 Curves, Jacobians, and abelian varieties (Amherst, MA, 1990),  55--125, Contemp. Math., 136, Amer. Math. Soc., 1992. Available online at
               \href{http://arxiv.org/alg-geom/9206008}{alg-geom/9206008}.
\bibitem[DL]{DL}  Donagi, R.; Livn\'e, R. {\em The arithmetic-geometric mean and isogenies for curves of higher genus}.
 Ann. Scuola Norm. Sup. Pisa Cl. Sci. (4)  28  (1999),  no. 2, 323-339.
Available online at \href{http://arxiv.org/alg-geom/9712027}{alg-geom/9712027}.
\bibitem[GH]{GH} Gross, B.; Harris, J. {\em Real algebraic curves}. Ann. Sci. \'Ecole Norm. Sup. (4) 14 (1981), no. 2, 157-182. Available online at
http://archive.numdam.org/article/ASENS\_1981\_4\_14\_2\_157\_0.pdf
\bibitem[Har]{Ha} Harris, J. {\em Galois groups of enumerative problems}. Duke Math. J., 46 (1979), no. 4, 685-724.
\bibitem[Har]{Har}  R. Hartshorne {\em Algebraic Geometry}
              Springer Verlag, GTM 52. (1977)
\bibitem[Hui]{Hui} Huisman, J. {\em  A group law on smooth real quartics having at least 3 real branches}.
J. Thor. Nombres Bordeaux 14 (2002), 249-256.
Available online at http://fraise.univ-brest.fr/\~{}huisman/recherche/publications/jq.html
\bibitem[Hum]{Hum} Humbert, G. {\em Sur la transformation ordinaire des fonctions abeliennes}.
              J. de math. (5) 7 (1901).
\bibitem[Jor]{Jo} Jordan, M.C. {\em Trait\'e des substitutions et des \'equations alg\'ebriques}. Gauthier-Villars, Paris, 1870.
\bibitem[Kon]{Ko} Konigsberger, L. {\em Uber die transformation der Abelschen Functionen erster ordnung}.
              J. reine angew. Math. 64 (1865) 17-42.
\bibitem[Leh]{Le} Lehavi, D. {\em A smooth plane quartic can be reconstructed from its bitangents}. Isr. J. Math. 146, 371-379 (2005).
            Available online at \href{http://arxiv.org/math.AG/0111017}{math.AG/0111017}.
\bibitem[LL]{LL} R. Lercier \& D. Lubicz : {\em A quasi quadratic time
algorithm for hyperelliptic curve point counting}. Available online
at
http://www.medicis.polytechnique.fr/~lercier/preprints/riemann.pdf.
\bibitem[LR]{LR} Lehavi, D; Ritzenthaler, C. {\em A proof of the arithmetic geometric mean formula in genus $3$ - a computer program}. Availbale online at \href{http://www.math.ohio-state.edu/~dlehavi}{http://www.math.ohio-state.edu/\~{}dlehavi}
\bibitem[Mes]{Me} Mestre, J.-F. {\em Lettre adress\'ee \`a Gaudry et Harley, D\'ecembre 2000}.
  Available online at \href{http://www.math.jussieu.fr/~mestre/}{http://www.math.jussieu.fr/\~{}mestre/}
\bibitem[Mum]{Mu} Mumford, D. {\em Prym varieties I}.
 Contributions to analysis (a collection of papers dedicated to Lipman Bers),  pp. 325--350.
Academic Press, 1974.
\bibitem[Pan]{Pa} S. Pantazis {\em Prym Varieties and the Geodesic Flow on $SO(n)$}
             Math. Ann. 273 p. 297-315 (1986).
\bibitem[Rec]{Re} Recillas, S. {\em Jacobians of curves with $g^1_4$'s are the Pryms of trigonal curves}.  Bol. Soc. Mat. Mexicana (2)  19  (1974), no. 1, 9-13.
\bibitem[Ric]{Ric} Richelot, F. {\em De transformatione integralium Abelianorum primi ordinis comentatio}.
              J. reine angew. Math. 16 (1837) 221-341.
\bibitem[Rit]{Rit} Ritzenthaler, C. {\em Probl\`emes arithm\'tiques relatifs \`a certaines familles de courbes sur les corps  nis}.
Th\`ese de Doctorat, Universit\'e Paris 7. Available online at \href{http://www.math.jussieu.fr/~ritzenth/}{http://www.math.jussieu.fr/\~{}ritzenth/}.
\bibitem[Sal]{salmon} G. Salmon : \textit{A treatise on the higher plane
    curves}, troisi{\`e}me {\'e}dition, Chelsea, (1879).
\end{thebibliography}
\end{document}